\documentclass{amsart}
\usepackage{amsmath, amsthm, amssymb, calrsfs, wasysym, verbatim, bbm, color, graphics, graphicx, geometry,enumerate,cite}
\usepackage{scrextend}
\usepackage[all,cmtip]{xy}
\usepackage{tikz}
\usepackage{tikz-cd}
\usetikzlibrary{patterns} 
\usepackage{hyperref}
\hypersetup{pdfborder = {0 0 0}}
\usepackage[nameinlink]{cleveref}
\usepackage{hyperref}
\hypersetup{urlcolor=blue, linktocpage, hyperindex=true, colorlinks=true,
	linkcolor=blue, citecolor=blue} 

\newtheorem{theorem}{Theorem}[section]
\newtheorem{lemma}[theorem]{Lemma}

\theoremstyle{definition}
\newtheorem{definition}[theorem]{Definition}
\newtheorem{proposition}[theorem]{Proposition}

\newtheorem{example}[theorem]{Example}
\newtheorem{corollary}[theorem]{Corollary}
\newtheorem{question}[theorem]{Question}

\theoremstyle{remark}
\newtheorem{remark}[theorem]{Remark}
\numberwithin{equation}{section}

 \newcommand{\ZZ}{\mathbb{Z}} 
\newcommand{\QQ}{\mathbb{Q}}
\newcommand{\CC}{\mathbb{C}}
\newcommand{\NN}{\mathbb{N}}

\newcommand{\ot}{\otimes}
\newcommand{\Ker}{\mathrm{Ker}}

\newcommand{\cV}{\mathcal{V}}

\newcommand{\cO}{\mathcal{O}}

\newcommand{\Aut}{\mathrm{Aut}}
\newcommand{\GL}{\mathrm{GL}}

\newcommand{\Pic}{\mathrm{Pic}}

\newcommand{\Vol}{\mathrm{Vol}}
\newcommand{\curProd}{C_1\times C_2\times C_3}

\newcommand\co{\colon\thinspace}

\crefname{claim}{Cliam}{Claim}
\crefname{corollary}{Corollary}{Conjecture}
\crefname{theorem}{Theorem}{Theorem}
\crefname{conjecture}{Conjecture}{Conjecture}
\crefname{definition}{Definition}{Definition}
\crefname{equation}{Equation}{Eq.}
\crefname{example}{Example}{Example}
\crefname{lemma}{Lemma}{Lemma}
\crefname{proposition}{Proposition}{Proposition}
\crefname{remark}{Remark}{Remark}
\crefname{equstion}{Remark}{Question}

\begin{document}
	
	\title{Automorphisms of threefolds of general type acting trivially in cohomology}
	
	\author{Hang Zhao}
	\address{School of Mathematical Sciences \\ 
		Peking University   \\ 
		Beijing\\
		China}
	\email{math\_hang@pku.edu.cn}
	

	
	\date{\today}
	
	\subjclass[2010]{14J50, 14J30}
	\keywords{Threefolds isogenous to a product, Numerically trivial automorphism}
	
	
	\maketitle
	\begin{abstract}
		Let $X$ be a minimal projective threefold of general type over $\mathbb{C}$ with only Gorenstein quotient singularities, and let $\mathrm{Aut}_{\mathbb{Q}}(X)$ be the subgroup of automorphisms acting trivially on $H^*(X,\mathbb{Q})$. In this paper, we show that if $X$ is of maximal Albanese dimension, then $|\mathrm{Aut}_{\mathbb{Q}}(X)|\leq 6$. Moreover, if $X$ is nonsingular and $K_X$ is ample, then $|\mathrm{Aut}_{\mathbb{Q}}(X)|\leq 5$.
		
		Seeking for higher-dimensional examples of varieties with nontrivial $\mathrm{Aut}_{\mathbb{Q}}(X)$, we concern $d$-folds $X$ isogenous to an unmixed product of curves. If $d=3$, we show that $\mathrm{Aut}_{\mathbb{Q}}(X)$ is a $2$-elementray abelian group whose order is at most $4$ under some conditions on their minimal realizations. Moreover, each of the possible groups can be realized. If $d\geq 3$, we give a sufficient condition for $\mathrm{Aut}_{\mathbb{Q}}(X)$ being trivial.
		
		Curiously, there exist examples of projective threefolds $X$ with terminal singularities and maximal Albanese dimension whose $\mathrm{Aut}_{\mathbb{Q}}(X)$ can have an arbitrarily large order.
	\end{abstract}
	\tableofcontents
	\section{Introduction}
	\label{intro}
	Let $X$ be a complex manifold, and $\mathrm{Aut}(X)$ be its group of holomorphic automorphisms. Consider the action of $\mathrm{Aut}(X)$ on the cohomology $H^*(X,A)$ of $X$, where $A=\QQ,\ZZ$; this gives a representation $\rho_{X,A}\co \mathrm{Aut}(X)\rightarrow\mathrm{GL}(H^*(X,A))$ defined by $\rho_{X,A}(\sigma)(\omega)=(\sigma^{-1})^*\omega$ for $\sigma\in\mathrm{Aut}(X)$ and $\omega\in H^*(X,A)$. The interesting question is that is the representation $\rho_{X,A}$ faithful?
	
	We say  that $X$ is \emph{rationally cohomologically rigidified} (resp. \emph{cohomologically rigidified}) if $\rho_{X,\QQ}$ (resp. $\rho_{X,\ZZ}$) is faithful. The connected component of the identity $\Aut(X)^0\subset\Aut(X)$ acts trivially on the cohomology, and is therefore contained in the kernel of $\rho_{X,\QQ}$. In general, those automomorphisms acting trivially on $H^*(X,\QQ)$ are called \emph{numerically trivial} and they form a subgroup of the (full) automorphism group, to be denoted by $\Aut_{\QQ}(X)$ in this paper. Thus, the group $\mathrm{Aut}(X)$ splits into two basic parts: its neutral component $\mathrm{Aut}(X)^0$, and its discrete image $\Aut(X)^*\subset\mathrm{GL}(H^*(X,\QQ))$. The group of connected coponents $\Aut(X)/\Aut(X)^0$ is an extension of $\Aut(X)^*$ by $\Aut_{\QQ}(X)/\Aut(X)^0$. It is interesting to study the structure of the group $\mathrm{Aut}_{\mathbb{Q}}(X)$. In particular, when $X$ is of general type, $\mathrm{Aut}_{\mathbb{Q}}(X)$ is a finite group.
	
	Many authors have recently studied the numerically trivial automorophism group of surfaces whose Kodaira dimension ranging from 0 to 2.
	\begin{itemize}
		\item for K3 surfaces \cite{BHPV04,BR75};
		\item for Enriques surfaces \cite{Muk10,MN84};
		\item for properly elliptic surfaces \cite{cai2009automorphisms,peters1979automorphisms};
		\item for surfaces of general type \cite{cai2004automorphisms,cai2006automorphisms,cai2007classification,peters1979holomorphic,peters1979automorphisms,Cai10,cai2012automorphisms,cai2018automorphisms,cai2013automorphisms}.
	\end{itemize} 
	For surfaces of general type, it turns out that nontrivial $\Aut_{\QQ}(S)$ occurs only for those with
	irregularity $q(S)\leq 2$ due to Cai, Liu, and Zhang, who prove the following theorem.
	
	\begin{theorem}[\cite{cai2013automorphisms}]\label{thm: surfaces AutQ}
		Let $S$ be a minimal surface of general type. Then we have the following results:
		\begin{enumerate}
			\item if $q(S)\geq 3$, then $S$ is rationally cohomologically rigidified;
			\item if $q(S)=2$, then $|\mathrm{Aut}_{\mathbb{Q}}(S)|\leq 2$, and the equality holds only if $S$ is a surface isogenous to an unmixed product of curves.
		\end{enumerate}
	\end{theorem}
	
	This paper aims to study the numerically trivial automorphism group of threefolds of general type with maximal Albanese dimension. First, we prove the following result.
	\begin{theorem}[Theorem \ref{thm: bound}]\label{thm: bound intro}
		Let $X$ be a minimal projective threefold of general type with only Gorenstein quotient singularities, assume that it is of maximal Albanese dimension. Then $|\mathrm{Aut}_{\mathbb{Q}}(X)|\leq 6$. Moreover, if $X$ is smooth and $K_X$ is ample, $|\mathrm{Aut}_{\mathbb{Q}}(X)|\leq 5$.
	\end{theorem}
	
	Our proof is inspired by \cite{cai2018automorphisms,cai2013automorphisms}. Set $\bar{X}:=X/\mathrm{Aut}_{\mathbb{Q}}(X)$. One can show that  the Albanese map $a_X\co X\rightarrow A_X$ of $X$ factors through the quotient map $\pi\co X\rightarrow \bar{X}$ and that $\chi(\omega_X)=\chi(\omega_{\bar{X}})$, see Lemma \ref{lem: quotient of X by AutQ}, (1). Since $X$ is a Gorenstein minimal threefold, the Bogomolov-Miyaoka-Yau inequality \cite[Theorem 1.1]{miyaoka1987chern} implies that
	$$\Vol(K_X)\leq 72\chi(\omega_X).$$ 
	Let $\bar{Y}$ be a suitable desingularization of $\bar{X}$, we can show that $\bar{Y}$ is of general type and of maximal Albanese dimension, see (2) of Lemma \ref{lem: quotient of X by AutQ}. By the generalized Severi inequality \cite{Bar15,ZT14}, we have 
	$$12\chi(\omega_{\bar{Y}})\leq\Vol(K_{\bar{Y}}).$$ Comparing volumes $\Vol(K_X)$ and $\Vol(K_{\bar{Y}})$, we obtain 
	$$|\mathrm{Aut}_{\mathbb{Q}}(X)|\Vol(K_{\bar{Y}})\leq \Vol(K_X).$$
	Combining the inequalities above, we get $|\mathrm{Aut}_{\mathbb{Q}}(X)|\leq 6$. The assumption that $X$ has Gorenstein singularities is necessary; we give a counterexample when $X$ has terminal singularities of Cartier index $2$  by constructing a series of threefolds $\{X_n\}_{n\in\mathbb{N}}$ such that $|\mathrm{Aut}_{\mathbb{Q}}(X_n)|$ can be arbitrarily large, see Example \ref{ex:Z2n}.
	
	It is well-known that compact Riemann surfaces of genus $g\geq2$ are rationally cohomologically rigidified. In Section \ref{sec:5.1}, we generalize this fact to varieties isogenous to an unmixed product of curves, introduced by Catanese in \cite{Cat00}, which is a quotient of a product of curves of genus at least $2$ by a finite group acting freely and diagonally.
	\begin{theorem}[Theorem \ref{thm: d-folds isogenous AutQ}]\label{thm: d-folds isogenous AutQ intro}
		Let $X$ be $d$-fold isogenous to an unmixed product of curves with $d\geq 3$, and let $(C_1\times \dots\times C_d)/G$ be its minimal realization. Suppose $g(C_i/G)\geq 1$ for all $1\leq i\leq d$. Set $K_i=\mathrm{Ker}(G\rightarrow\mathrm{Aut}(C_i))$. If there is a pair $(i,j)$ with $j\neq i$ such that $g(C_i/G)\geq 2$ and $g(C_j/K_i)\geq 2$, then $\mathrm{Aut}_{\mathbb{Q}}(X)$ is trivial. 
	\end{theorem}
	
	In general, we first consider the case that a nonsingular projective variety $Y$ of dimsion $d\geq 3$ which admits a higher irrational pencil $g\co Y\to D$ where $D$ is a smooth curve. Let $\sigma$ be a nontrivial automorphism of $Y$ such that $g\circ\sigma=g$. If $\sigma$ induces a trivial action on $H^0(Y,\omega_Y)$, then its restriction ot $F$ induces the identity on $H^0(F,\omega_F)$, where $F$ is a general fibre of $g$. Let $\sigma_F$ be the resctriction of $\sigma$ on $F$. We can use lower dimension result on pair $(F,\sigma_F)$ to get a higher dimension result on pair $(Y,\sigma)$, in fact, we can show that $o(\sigma)\leq o(\sigma_F)$, see Lemma \ref{lem: g geq 1 AutQ to Aut F}. This result is a generalization of \cite[Lemma 2.1]{Cai12}. Then we can use a induction procedure to the case that there are a sequence of higher irrational pencils $g_j\co F_j\to D_j$ such that each $F_{j+1}$ is a general fibre of $g_j$ for $0\leq j\leq \dim Y-2$ where $F_0=Y$, and show that there is no nontrivial automorphism of $F_{d-2}$ which is a successive restriction of an automorphism $\tilde{\sigma}_0$ of $Y$ acting trivially on $H^0(Y,\omega_Y)$, see Corollary \ref{cor: Situation *}. In particular, applying this result to the case that $X$ is as in Theorem \ref{thm: d-folds isogenous AutQ intro} shows the theoerm.
	
	It is worth pointing out that Theorem \ref{thm: d-folds isogenous AutQ intro} is not valid for the case that $q(X)\geq d+1$, which is different from that of irregular surfaces of general type \cite[Theorem 1.4]{cai2013automorphisms}. In these papers \cite{cai2013automorphisms,liu2018surfaces}, both authors construct a series of surfaces $X$ of general type isogenous to a product with $q(X)=2$ such that $\mathrm{Aut}_{\mathbb{Q}}(X)\cong\ZZ_2$. In their example, $X$ is a quotient of $C\times D$ by a finite group $G$, where $C$ and $D$ are curves with faithful group actions of $G$. 
	
	Seeking for higher-dimensional examples of varieties with nontrivial $\mathrm{Aut}_{\mathbb{Q}}(X)$, we consider threefolds  $X$ isogenous to a product of curves. A new phenomenon occurs in this case: let $(C_1\times C_2\times C_3)/G$ be the minimal realization of $X$, then the group $G$ probably does not act faithfully on each curve $C_i$ for $i=1,2,3$. Denote the subgroup of $G$ acting trivially on $C_i$ by $K_i$; the appearance of $K_i$ is a difficulty for studying the structure of $\mathrm{Aut}_{\mathbb{Q}}(X)$. Suppose $G$ is abelian, and all $K_i$ are cyclic groups, then we can show that $\mathrm{Aut}_{\mathbb{Q}}(X)$ is a $2$-elementary abelian group. Concretely, we have the following theorem.
	\begin{theorem}[Theorem \ref{thm: threefold isogenous AutQ}]\label{thm: threefold isogenous AutQ intro}
		Let $X$ be threefold isogenous to an unmixed product of curves, and let $(C_1\times C_2\times C_3)/G$ be its minimal realization. Suppose $g(C_i/G)\geq 1$ for all $1\leq i\leq 3$. Then we have
		\begin{enumerate}[(1)] 
			\item If there is a pair $(i,j)$ with $j\neq i$ such that $g(C_i/G)\geq 2$ and $g(C_j/K_i)\geq 2$, then $\mathrm{Aut}_{\mathbb{Q}}(X)$ is trivial;
			\item if for any $1\leq i\leq 3$ with $g(C_i/G)\geq 2$, we have $g(C_j/K_i)=1$ for all $j\neq i$, then $\mathrm{Aut}_{\mathbb{Q}}(X)\cong(\mathbb{Z}_2)^k$ with $k=0,1$;
			\item if for all $1\leq i\leq 3$, we have $g(C_i/G)=1$, and suppose that the group $G$ is an abelian group, and $K_i$ is a cyclic group for all $1\leq i\leq 3$, then we have $\mathrm{Aut}_{\mathbb{Q}}(X)\cong (\mathbb{Z}_2)^k$ with $k=0,1,2$.
		\end{enumerate} 
	\end{theorem}
	
	We give only the main ideas of the proof. Theorem \ref{thm: threefold isogenous AutQ intro}, (1) follows directly from Theorem \ref{thm: d-folds isogenous AutQ intro}. We apply the induction procedure mentioned above to derive Theorem \ref{thm: threefold isogenous AutQ intro}, (2), see Corollary \ref{cor: class II}. To prove Theorem \ref{thm: threefold isogenous AutQ intro}, (3), we first show that $\mathrm{Aut}_{\mathbb{Q}}(X)$ is determined by the algebraic data associated with the minimal realization of $X$, see Definition \ref{def: algebraic data} for the definition of algebraic data and Lemma \ref{lem: representation of LfX}. We next show that the group $\mathrm{Aut}_{\mathbb{Q}}(X)$ can be embedded into an abstract $2$-elementary abelian group, see Theorem \ref{thm: G/KD is 2 elementary abelian group}. Finally, our assertion follows from the bound of $|\mathrm{Aut}_{\mathbb{Q}}(X)|$ in Theorem \ref{thm: bound intro}.
	
	In Section \ref{sec:6}, we construct some examples of threefolds $X$ isogenous to an unmixed product of curves with $\mathrm{Aut}_{\mathbb{Q}}(X)\cong \ZZ_2$ and $\ZZ_2\times \ZZ_2$. 
	
	\section{Notations and preliminaries}
	\label{sec:1}
	We work over the complex numbers throughout the paper.
	Varieties are always assumed to be normal and quasi-projective; a threefold $X$ is a projective variety of dimension $3$.
	
	Let $\Aut(X)$ be the holomorphic automorphism group of $X$. Let $f\co X\rightarrow Y$ be a surjective morphism to a variety $Y$ with connected fibres, we set 
	\begin{description}
		\item[$\mathrm{Aut}(X/Y)=\{\sigma\in\mathrm{Aut}(X)|f\circ\sigma=f\}$,] the relative automorphism group over $Y$.
	\end{description}
	
	For a finite group $G$, we denote
	\begin{description}
		\item[$G^*$] the set of irreducible characters on $G$.
	\end{description} When $G$ is abelian, $G^*$ is isomorphic to $G$, and is called the dual group of $G$.  For an element $g$ of $G$, we denote the order of $g$ by $o(g)$. 
	
	For a variety $X$ with a faithful group action of $G$ and a cohomology group $H$ of $X$, we set 
	\begin{description}
		\item[$X^g=\{x\in X|g(x)=x\},$] the fixed loci of an element $g\in G$;
		\item[$H^G=\{\omega\in H| g^*\omega=\omega\},$] the $G$-invariant subgroup of $H$.
	\end{description} 
	
	\subsection{Gorenstein quotient singularity}
	\label{sec:1.1}
	For a normal variety $X$ of dimension $n$ with singularities, we denote its nonsingular part by $X_0$, and then the canonical sheaf $\omega_{X_0}$ make sense. Let $j\co X_0\rightarrow X$ be the inclusion morphism, and we define the sheaves $\tilde{\Omega}^i_X=j_*(\Omega^i_{X_0})$ of $X$ for $1\leq i\leq n$. The \emph{canonical divisor} $K_X$ on $X$ is a Weil divisor such that $$\mathcal{O}_X(K_X)=\omega_{X}:=\tilde{\Omega}^n_X,$$ which is well-defined up to linear equivalence. We say that $X$ has \emph{Gorenstein singularity} if its canonical divisor $K_X$ is a Cartier divisor. For the definitions of rational singularity, we recommend references \cite[Section 2]{Ko13} and \cite[Section 6.2]{Is14}. 
	\begin{remark}\label{rem: rational singularity} 
		If the variety $X$ has rational singularities, for example, quotient singularities, for any nonsingular resolution $f\co  Y\rightarrow X$, we have $R^if_*\cO_Y=0$ for $i>0$ and $f_*\cO_Y=\cO_X$. It follows that $\chi(X,\cO_X)=\chi(X,f_*\cO_Y)=\chi(Y,\cO_Y)$.  Since rational singularities are Cohen-Macaulay (see \cite[Theorem 6.2.14]{Is14}), using the Serre duality, we get $\chi(X,\omega_X)=\chi(Y,\omega_Y)$.
	\end{remark}
	
	\subsection{Volume of divisors}\label{sec:1.2}
	\begin{definition}\label{def: volume}
		Let $X$ be a projective variety of dimension $n$, and let $D$ be an integral divisor on $X$. The volume of $D$ is defined to be the non-negative real number $$\Vol_X(D)=\limsup_{m\rightarrow\infty}\frac{n!h^0(X,\cO_X(mD))}{m^n}$$
		If $D$ is a $\QQ$-divisor, the volume of $D$ is defined as $\Vol_X(D)=\frac{1}{a^n}\Vol_X(aD)$ for some $a\in\NN$ such that $aD$ is integral.
	\end{definition}
	
	\begin{proposition}\cite[Proposition 2.2.43]{LR04}\label{prop: volume Vol K_X'= Vol K_X}  
		Let $X$ be a normal projective variety of dimension $n$.
		If $X$ has canonical singularities and $\nu\co X'\rightarrow X$ is a nonsingular resolution of $X$, then $$\Vol_{X'}(K_{X'})=\Vol_X(K_X)$$
	\end{proposition}
	\begin{remark}\label{rem: Vol D leq Vol D+F}
		Note that $\Vol_X(D)>0$ if and only if $D$ is big. If $D$ is nef, then it follows from the asymptotic Riemann-Roch that $\Vol_X(D)=D^n$. If $F$ is an effective divisor on $X$, we have $\Vol_X(D)\leq\Vol_X(D+F)$.\end{remark}
	
	\subsection{Galois covers of curves}\label{sec:1.3}
	We recall some facts about the Galois covers of algebraic curves, and refer the reader to \cite[Section 2]{Br90}, \cite[Chapter 3]{Bre00} and \cite[Section 1]{Pol08} for more details.
	\begin{definition}\label{def: Riemann's existence theorem}
		Let $G$ be a finite group and let $$g'\geq 0,\quad m_r\geq m_{r-1}\geq\dots\geq m_1\geq2$$ be integers. A generating vector for $G$ of type $[g';m_1,\dots,m_r]$ is
		a sequence of elements $$\cV:=(\eta_1,\dots,\eta_{2g'};\sigma_1,\dots,\sigma_r )$$ such that the following conditions are satisfied,
		\begin{enumerate}[(1)]
			\item $G$ is generated by the entries of the sequence $\cV$;
			\item $o(\sigma_i)=m_i$ for all $1\leq i\leq r$;
			\item $\sigma_1 \cdots\sigma_r\cdot\prod_{i=1}^{g'}[\eta_i ,\eta_{g'+i} ]=1_{G}$;	
		\end{enumerate} Moreover, if the following Riemann-Hurwitz relation holds:$$2g-2=|G|(2g'-2+\sum_{i=1}^r(1-\frac{1}{m_i})),$$
		then there exists a compact Riemann surface $C$ of genus $g$ and a $G$-cover $C\rightarrow D$ branched in $r$ points $P_1,\dots,P_r$ with ramification indexes $m_1,\dots,m_r$ respectively, where $g'$ is the genus of $D$. The subgroups $g\langle \sigma_i\rangle g^{-1}$ provide all nontrivial stabilizers of the $G$-action on $C$ for $g\in G$. Set $\Sigma:=\bigcup_{i=1}^r\bigcup_{g\in G}g\langle\sigma_i\rangle g^{-1}$. \end{definition}
	
	In the terminology of representations, let $\varphi\co G\rightarrow\GL(H^{1,0}(C))$ be the representation given by the action of $G$ on $C$. For any irreducible representation $\rho\co G\rightarrow\GL(V)$ with $V$, a finite-dimensional vector space over $\CC$. We denote by $N_{i,k}(\rho)$ the multiplicity of $\mu_{m_i}^k$ as an eigenvalue of $\rho(\sigma_i)$, where $\mu_{m_i}:=\mathrm{exp}(\frac{2\pi\sqrt{-1}}{m_i})$ and $0\leq k\leq m_i-1$. For the character $\chi_{\varphi}$ afforded from $\varphi$, we have the following formula.
	\begin{theorem}[Chevalley-Weil formula \cite{chevalley1934verhalten,gleissner2016threefolds}]\label{thm: Chevalley-Weil formula} 
		In the notations above, let $1_G$ be the trivial character of degree $1$ on $G$. For each $\chi\in G^*$ and $\rho$ is its corresponding representation, it holds $$\langle\chi,\chi_{\varphi}\rangle = \chi(1_G)(g'-1) + \sum_{i=1}^r\sum_{k=1}^{m_i-1}\frac{k\cdot N_{i,k}(\rho)}{m_i}+\langle\chi,1_G\rangle.$$ \end{theorem}
	
	\begin{remark}\label{rem: condition for H^1chi neq 0} 
		Under the assumption above, suppose $G$ is abelian and $g'\geq1$. Then for each $\chi\in G^*\setminus\{1_G\}$, $\dim_{\CC}H^{1,0}(C)^{\chi}>0$ if and only if the following holds \begin{equation*}\label{eq2} \sigma\notin\Ker(\chi)\text{ for some }\sigma\in\Sigma, \text{i.e., }\sigma=\sigma^k_i \text{ for some } 1\leq i\leq r,1\leq k\leq m_i-1.\end{equation*}
	\end{remark}
	
	\subsection{Threefolds isogenous to a product of curves}\label{sec:1.4}
	According to Catanese \cite{Cat00}, a vareity $X$ is said to isogenous to a higher porduct of curves if it admits a finite unramified covering which is isomorphic to  a product of curves of genus $\geq 2$. There is a product of curves $C_1\times\dots\times C_d$ together with finite group $G$ such that $C$ acts freely on $C_1\times\dots \times C_d$. The pair $(C_1\times\dots\times C_d, G)$ or the quotient $(C_1\times\dots\times C_d)/G$ is called a \emph{minimal realization} of $X$ if $K_i\cap K_j=\{1\}$ for all $1\leq i<j\leq d$, where $K_i$ is the normal subgroup of $G$ which acts trivially on $C_i$. $X$ is said to be of \emph{unmixed type} if $G$ acts diagonally on $C_1\times\dots\times C_d$, and is also said to \emph{isogenous to an unmixed product of curves}.
	
	\begin{definition}\label{def: algebraic data}
		Let $G$ be a finite abelian group with $K_1,K_2,K_3$ three normal subgroups, and let $\cV_i$ be a generating vector for the quotient group $G/K_i$ for $1\leq i\leq3$. The 7-tuple $\mathcal{A}=(G,K_1,K_2,K_3,\cV_1,\cV_2,\cV_3)$ is called an algebraic datum for $G$ if the following conditions are satisfied:
		\begin{enumerate}[(1)]
			\item $K_i\cap K_j=\{1_G\}$ for each $1\leq i<j\leq3$,
			\item $\Sigma_1\cap\Sigma_2\cap\Sigma_3=\{1_G\}$.
		\end{enumerate}
		Where $\Sigma_i$ is the union of nontrivial stabilizers of the $G$-action on each factor $C_i$ for $1\leq i\leq 3$.
	\end{definition}
	
	\begin{remark}\label{rem: construction of threefold isogenous}
		A threefold isogenous to an unmixed product of curves is determined by an algebraic datum $\mathcal{A}$ described above.  For each algebraic data $(G/K_i,\cV_i)$, by the Riemann existence theorem, there exists an algebraic curve $C_i$ with a faithful group action of $G/K_i$. We denote by $\psi_i\co G\rightarrow \Aut(C_i)$, the action of $G$ on $C_i$; and we have $K_i=\Ker(\psi_i)$. The homomorphisms $\psi_1,\psi_2$ and $\psi_3$ induce a $G$-action on the product $\curProd$:
		$$g(x_1,x_2,x_3)=(\psi_1(g)x_1,\psi_2(g)x_2,\psi_3(g)x_3),$$
		where $g\in G$ and $(x_1,x_2,x_3)\in C_1\times C_2\times C_3$. The second condition, which is called the freeness condition, ensures that the action of $G$ on $C_1\times C_2\times C_3$ is free. Therefore, the quotient $X=(C_1\times C_2\times C_3)/G$ is a threefold isogenous to a product of curves.
	\end{remark}
	
	\subsection{Characters of finite abelian groups} \label{sec:1.5}
	For basic definitions of the representation theory of groups, we refer to the books \cite{IS76,Se77}.
	
	Let $G$ be a finite abelian group with the identity $1$, and let $H$ be a subgroup of $G$. The restriction map $G^*\rightarrow H^*$ is a surjective group homomorphism, for $\chi\in G^*$, we denote the restriction of $\chi$ on $H$  by $\chi_H$. Since $(G/H)^*\cong\Ker(G^*\rightarrow H^*)$, we may identify $(G/H)^*$ with the subset of characters of $G$ whose restriction on $H$ is trivial. If $G$ is a cyclic group with a generator $e$, we say that a character $\chi$ of $G$ is \emph{primitive} if $\chi(e)$ is a $o(e)$-th primitive root of unit, so a primitive character $\chi$ is a generator of the dual group $G^*$. For a primitive character $\chi$ of $G$ and $g\in G$, we have $\chi(g)\neq 1$ iff $g\neq1$. Moreover, $\chi(g)$ is a $o(g)$-th primitive root of unit. 
	
	Let $\rho\co G\to \GL(V)$ be a linear representation of $G$ over $\CC$, and let $\chi$ be the character of the representation $\rho$ given by $\chi(s)=\mathrm{Tr}(\rho(s))$ for each $s\in G$. Then $V$ decomposes into a direct sum of irreducible representations:
	$$V=V^{\chi_1}\oplus \dots\oplus V^{\chi_k},$$
	where $V^{\chi_i}$ is sum of irrducible representation with character $\chi_i$ for all $1\leq i\leq k$. Set $n_i=\dim V^{\chi_i}$, then we can write $\chi=\sum_{i=1}^kn_i\chi_i$.
	
	Let $G=H\times K$ be a direct product of finite groups and let $\varphi$ and $\theta$ be characters on $H$ and $K$, respectively. We define a character $\chi=\varphi\times\theta$ of $G$ by $\chi(hk)=\varphi(h)\theta(k)$ for $h\in H$ and $k\in K$. Since we have $H\cong G/K$, there is a corresponding character $\hat{\varphi}$ of $G$ such that $K\subseteq\Ker(\hat{\varphi})$ and $\hat{\varphi}(hk)=\varphi(h)$. Similarly, there is a corresponding character $\hat{\theta}$ of $G$ such that $H\subseteq\Ker(\hat{\theta})$ and $\hat{\theta}(hk)=\theta(k)$. It follows that $\varphi\times\theta=\hat{\varphi}\hat{\theta}$. Moreover, the characters $\varphi\times\theta$ for which $\varphi$ and $\theta$ being irreducible are exactly the irreducible characters of $G$. Let $V$ and $W$ be linear representations of $H$ and $K$, respectively. Let $V=\oplus_{1\leq i\leq k}V^{\chi_i}$ and $W=\oplus_{1\leq j\leq l}W^{\psi_j}$ be corresponding decompositions. Then $V\otimes W$ is a linear representation of $G$ with decomposition:
	\begin{align}\label{eq: decomposition of V otimes W}
		V\otimes W&=\bigoplus_{1\leq i\leq k}\bigoplus_{1\leq j\leq l}V^{\chi_i}\otimes W^{\psi_j} \nonumber\\
		&=\bigoplus_{i,j}(V\otimes W)^{\chi_i\times\psi_j} 
	\end{align}
	
	Let $G$ be a finite group, not necessarily abelian, let $H\subseteq G$ be a subgroup, and let $\varphi$ be a character of $H$. We define the \emph{induced character} $\varphi^G$ of $G$ by
	\begin{equation}\label{eq: defn of induced character}
		\varphi^G(g)=\frac{1}{|H|}\sum_{x\in G}\varphi^0(xgx^{-1}),
	\end{equation}
	where $\varphi^0$ is defined by $\varphi^0(h)=\varphi(h)$ if $h\in H$ and $\varphi^0(y)=0$ if $y\notin H$. If $G$ is abelian and $g\notin H$, then we have $\varphi^G(g)=0$. On the other hand, we may write $\varphi^G=\sum_{i=1}^k n_i\chi_i$ with $\chi_i$ is an irreducible character of $G$, we call each $\chi_i$ a \emph{constituent} of $\varphi^G$.  
	
	Next, we present some technical results required in the proof of Theorem \ref{thm: threefold isogenous AutQ}.
	
	\begin{proposition}\label{prop: consituent of induced character such that c chi neq 1}
		Let $G$ be a finite abelian group, $H$ a proper subgroup of $G$, and let $\varphi$ a character of $H$. If $g\notin H$, then for any root of unit $c$ there is a constituent $\chi$ of $\varphi^G$ such that $c\chi(g)\neq 1$.
	\end{proposition}
	
	\begin{proof} Since $g\notin H$, by formula (\ref{eq: defn of induced character}), we have $\varphi^G(g)=0$. Write $\varphi^G=\sum_{i=1}^k n_i\chi_i$. Suppose $c\chi_i(g)=1$ for all $\chi_i$ constituent of $\varphi^G$, then we have $$\varphi^G(g)=\sum_{i=1}^k n_i\chi_i(g)=\bar{c}[G:H]\neq0,$$ which is a contradiction. 
		\qed\end{proof}
	
	\begin{lemma}\label{lem: construction of bases characters} 
		Let $K_1$ and $K_2$ be two cyclic groups, let $H=K_1\times K_2$ be their product, and let $K_3$ be a cyclic subgroup of $H$ such that $K_1\cap K_3=\{1\}$ and $K_2\cap K_3=\{1\}$. Then there are nontrivial characters $\alpha_1,\alpha_2$ of $H$ satisfying $$K_1\subseteq\Ker(\alpha_2),K_2\subseteq\Ker(\alpha_1)\text{ and }K_3\subseteq\Ker(\alpha_1\alpha_2).$$ In particular, we can choose $\alpha_i$ such that its restriction on $K_i$ is a primitive character for each $i=1,2$. 
	\end{lemma}
	\begin{proof} Choose a generator $k_i$ of sugroup $K_i$ for each $i=1,2,3$. We can write $k_3=k_1^ak_2^b$ for some integers $a$ and $b$. Set $$n_1=o(k_1),n_2=o(k_2),m_1=o(k_1^a)\text{ and }m_2=o(k_2^b).$$ By assumption $K_3\cap K_1=K_3\cap K_2=\{1\}$, we conclude that $o(k_3)=m_1=m_2$, and set $m=o(k_3)$. We can choose two characters $\beta_1,\beta_2$ of $H$ satisfying $$\beta_i(k_j)=e^{2\delta_{ij}\pi\sqrt{-1}/n_j}$$ for each $1\leq i< j\leq 2$. Choose integers $s_1,s_2$ and set $\alpha_1=\beta_1^{s_1},\alpha_2=\beta_2^{s_2}$, to require $K_3\subseteq\Ker(\alpha_1\alpha_2)$, we note that $\alpha_1(k_3)\alpha_2(k_3)=e^{2(s_1+s_2)\pi\sqrt{-1}/m}=1$. Therefore it is sufficient to take $\alpha_1=\beta_1^s$ and $\alpha_2=\beta_2^{-s}$ for some integer $s$ that is relatively prime to $m$. To require $(\alpha_i)_{K_i}$ being primitive, it is sufficient to take $s=1$, which is due to the fact that the restriction of $\beta_i$ on $K_i$ is primitive. It is clear that $K_1\subseteq\Ker(\alpha_2)$ and $K_2\subseteq\Ker(\alpha_1)$, which completes the proof.
		\qed\end{proof}
	
	\begin{lemma}\label{lem: if max geq 3 construction of adm characters} 
		Let $K_1$ and $K_2$ be two cyclic groups, and let $H=K_1\times K_2$ be their product. Fix nontrivial element $g_i\in K_i$ for each $i=1,2$, and an element $h\in H$. Write $h=h_1h_2$ for $h_i\in K_i$. Let $\alpha_1,\alpha_2\in H^*$ be two characters satisfying the following conditions:
		\begin{enumerate}[(1)]
			\item $K_1\subseteq\Ker(\alpha_2),K_2\subseteq\Ker(\alpha_1)$;
			\item the restriciton of $\alpha_i$ on $K_i$ is primitive for each $i=1,2$.
		\end{enumerate}
		Set $I_i:=\{\chi\in H^*|\chi(g_i)\neq 1\}$ for each $i=1,2$ and $J:=\{\chi\in H^*|\chi(h)\neq1\}$.\\
		\indent	If $\max\{o(g_1),o(g_2),o(h_1),o(h_2)\}\geq 3$, then we have
		$$\begin{cases}
			\alpha_1^s\in I_1\text{ and }\alpha_1^s\alpha_2\in J & \text{ if }o(g_1)\geq 3\text{ or }o(h_1)\geq 3\\
			\alpha_2^s\in I_2\text{ and }\alpha_1\alpha_2^s\in J & \text{ if }o(g_2)\geq 3\text{ or }o(h_2)\geq 3
		\end{cases}$$
		for some $s=1,2,3,5$. 
	\end{lemma}
	\begin{proof}
		If $o(g_1)\geq 3$, then $\alpha_1^s\in I_1$ for $s=1,2$. Suppose $$\alpha_1(h)\alpha_2(h)=\alpha_1^2(h)\alpha_2(h)=1,$$ we have $\alpha_1(h_1)=\alpha_1(h)=\alpha_2(h)=\alpha_2(h_2)=1$, which contradicts to condition (2). Therefore, $\alpha_1^s\alpha_2\in J$ for some $s=1,2$. For the same reason, if $o(g_2)\geq 3$, we have $\alpha_1\alpha_2^s\in J$ for some $s=1,2$ with $\alpha_2^s\in I_2$. 
		
		Now we suppose that $o(g_1)=o(g_2)=2$ and $o(h_1)\geq 3$. If $o(h_1)=3$, we have $\alpha_1^s\alpha_2\in J$ and $\alpha_1^s\in I_1$ for some $s=1,5$. Otherwise, we have $$\alpha_1(h)\alpha_2(h)=\alpha_1^5(h)\alpha_2(h)=1,$$ and it follows that $\alpha_1^4(h)=1$. Note that  $\alpha_1^4(h)=\alpha_1(h_1)\neq 1$ by condition (2); we have a contradiction. If $o(h_1)\geq 4$, for the same reason, we have $\alpha_1^s\alpha_2\in J$ and $\alpha_1^s\in I_1$ for some $s=1,3$. Similarly, if $o(g_1)=o(g_2)=2$ and $o(h_2)\geq 3$, then $\alpha_1\alpha_2^s\in J$ and $\alpha_2^s\in I_2$ for some $s=1,3,5$. 
		\qed\end{proof}
	
	\begin{lemma}\label{lem: if o geq 3 construction of adm characters} 
		Let $G$ be a finite abelian group with the identity $1$, and let $g_1,g_2,h$ be three nontrivial elements in $G$. Set 
		\begin{equation*}
			J:=\{\phi\in G^*|\phi(g_1)\neq1\text{ and }\phi(g_2)\neq 1\}\text{ and }I:=\{\chi\in G^*|\chi(h)\neq1\}.
		\end{equation*} 
		If $o(h)\geq 3$, then we have $\phi\chi\in J$ for some $\phi\in J$ and $\chi\in I$.
	\end{lemma}
	\begin{proof} Let $H_i=\langle g_i\rangle$ for $i=1,2$ and $H=H_1H_2$. If $h\notin H$, then there is a character $\chi\in (G/H)^*$ such that $\chi(h)\neq 1$. Choose a character $\alpha\in H^*$ such that $\alpha(g_1)\neq 1$ and $\alpha(g_2)\neq 1$. Since the restriction map $G^*\rightarrow H^*$ is surjective, there is a character $\phi\in G^*$ such that $\phi_H=\alpha$.  Therefore, we have $\chi(g_i)\phi(g_i)=\phi(g_i)\neq 1$ for all $i=1,2$. Hence $\chi\in I$, $\phi\in J$ and $\phi\chi\in J$.
		
		Now we suppose that $h\in H$. The inclusion relation between subgroups $H_1,H_2$ and $H$ is one of the following forms:
		\begin{enumerate}
			\item $H_2\leq H_1$;
			\item $H_1$ and $H_2$ are distinct and their intersection $H_1\cap H_2$ is a proper subgroup in both $H_1$ and $H_2$.
		\end{enumerate}
		
		In the former case, we have $H=H_1$. Let $m$ be the order of $g_1$. Since $o(h)\geq 3$ and $h\in H$, we have $m\geq 3$. Let $\alpha\in H^*$ be the character such that $\alpha(g_1)$ is a $m$-th primitive root of unit. If $o(g_2)\neq 2$, we can take $\phi=\chi=\alpha$ such that $\phi\in J,\chi\in I$ and that $\phi\chi=\alpha^2\in J$; if $o(g_2)=2$, then $m\geq4$, we can take $\phi=\alpha,\chi=\alpha^2$ such that $\phi\in J,\chi\in I$ and that $\phi\chi=\alpha^3\in J$.
		
		In the latter case. Let $m_i$ be the order of $g_i$. We can write $h=g_1^ag_2^b$ for some integers $a$ and $b$. Since $o(h)\geq 3$, we have $\max\{o(g_1^a),o(g_2^b)\}\geq 3$. Now suppose $o(g_1^a)\geq 3$, so $m_1\geq 3$. If $H_1\cap H_2=\{1\}$, then there are characters $\alpha_i\in H^*$ such that $\alpha_i(g_j)=e^{\delta_{ij}2\pi \sqrt{-1}/m_j}$. We can take $\phi=\alpha_1\alpha_2$ and $\chi=\alpha_1$, it is easy to verify that $\phi\in J,\chi\in I$ and $\phi\chi\in J$. If $H_1\cap H_2\neq\{1\}$. Let $\alpha$ be a character of $H$ such that $\alpha(g_1)=e^{2\pi\sqrt{-1}/m_1}$, and let $n$ be the smallest integer such that $g_2^n\in H_1\cap H_2$, by assumption $n\geq2$. We can choose $\alpha$ such that $\alpha(g_2)^n=e^{2\pi\sqrt{-1}/d}$, where $d=m_2/n$, by assumption $d\geq 2$, so we have $m_2\geq 4$ and $\alpha(g_2)\neq 1$. Take $\phi=\chi=\alpha$, it is easy to verify that $\phi\in J,\chi\in I$ and $\phi\chi\in J$.
		\qed\end{proof}
	
	\section{The bound of $|\mathrm{Aut}_{\mathbb{Q}}(X)|$}\label{sec:3}
	\begin{theorem}\label{thm: bound} 
		Let $X$ be a minimal projective threefold of general type with only Gorenstein quotient singularities, assume that it is of maximal Albanese dimension. Then $|\mathrm{Aut}_{\mathbb{Q}}(X)|\leq 6$. Moreover, if $X$ is smooth and $K_X$ is ample, $|\mathrm{Aut}_{\mathbb{Q}}(X)|\leq 5$.
	\end{theorem}
	
	For the proof of the above theorem we need the following lemma.
	\begin{lemma}\label{lem: quotient of X by AutQ}  
		Let $X$ be a threefold as in Theorem \ref{thm: bound}, and set $G=\mathrm{Aut}_{\mathbb{Q}}(X)$. Then we have:
		\begin{enumerate}[(1)]
			\item the Albanese map $a_X\co X\rightarrow A_X$ of $X$ factors through the quotient map $\pi\co X\rightarrow \bar{X}$ and that  $\chi(\omega_X)=\chi(\omega_{\bar{X}})$;
			\item the quotient $X/G$ is of general type and of maximal Albanese dimension.
	\end{enumerate}\end{lemma}
	
	\begin{proof}To prove (1), let $\sigma\in G$ be a nontrivial automorphism. We claim that $X^{\sigma}\neq\varnothing$. Let $k$ be the smallest integer such that $X^{\sigma^k}\neq\varnothing$, then the quotient map $$\pi':X'=X/\langle\sigma^k\rangle\longrightarrow X''=X/\langle\sigma\rangle$$ is \'{e}tale. It follows that 
		\begin{equation}\label{eq: chi X'= ord chi X''}
			\chi(\omega_{X'})=|\langle\sigma\rangle/\langle\sigma^k\rangle|\chi(\omega_{X''})
		\end{equation} 
		Since $X$ has only quotient singularities, the spectral sequence
		$$E_1^{p,q}:=H^q(X,\tilde{\Omega}_X^p)\Rightarrow H^{p+q}(X,\mathbb{C})$$
		degenerates at $E_1$ page \cite[Theorem 2.43]{Pet08}. 
		And since $\sigma$ induces trivial action on $H^i(X,\CC)\cong H^i(X,\QQ)\otimes\CC$, it also induces identity action on $H^i(X,\omega_X)$ for all $i\geq 0$. It follows that
		\begin{align*}&H^i(X'',\omega_{X''})=H^i(X,\omega_X)^{\langle\sigma\rangle}=H^i(X,\omega_X)\\
			&H^i(X',\omega_{X'})=H^i(X,\omega_X)^{\langle\sigma^k\rangle}=H^i(X,\omega_X)
		\end{align*} for all $i\geq 0$, and that 
		\begin{equation}\label{eq: chi X= chi X' = chi X''}
			\chi(\omega_X)=\chi(\omega_{X'})=\chi(\omega_{X''}).
		\end{equation}
		Since $X$ is a Gorenstein minimal threefold of general type, by the Bogomolov-Miyaoka-Yau inequality \cite[Theorem 1.1]{miyaoka1987chern}, $0<K_X^3\leq 72\chi(\omega_X)$, we have $\chi(\omega_X)>0$. Combining this with the two equalities (\ref{eq: chi X'= ord chi X''})  and (\ref{eq: chi X= chi X' = chi X''}) we have $|\langle\sigma\rangle|=|\langle\sigma^k\rangle|$. Therefore, $X^{\sigma}\neq\varnothing$.
		
		Let $a_X\co X\rightarrow A_X$ be the Albanese map of $X$. Notice that $G$ induces identity on $H^1(X,\cO_X)$; for any $\sigma\in G$, the induced map $\bar{\sigma}$ on $A$ is a translation. Since $X^{\sigma}\neq\varnothing$ and $a(X^{\sigma})\subseteq A^{\bar{\sigma}}$, $\bar{\sigma}$ must be the identity map. Then the quotient map $X\longrightarrow \bar{X}$ factors through the Albanese map $a_X$ of $X$.
		
		By (1), there is a commutative diagram
		$$\xymatrix{
			X    \ar[rr]^{\pi} \ar[dr]_{a_X}  & &  \bar{X}  \ar[dl]^{a'}\\
			& A_X   &     }$$ According to the universal property of Albanese map $a_X$, the induced map $a'$ is the Albanese map of $\bar{X}$. As the map $a_X$ is generically finite onto its image and the quotient map $\pi$ is finite, we have that $\bar{X}$ is of maximal Albanese dimension. By the generic vanishing theorem (see \cite{green1987deformation,green1991higher}), for a general $\alpha\in\Pic^0(A)$, $h^i(\bar{X},\omega_{\bar{X}}\otimes a'^*\alpha)=0$ for all $i>0$, and hence $\chi(\omega_{\bar{X}})=h^0(\bar{X},\omega_{\bar{X}}\otimes a'^*\alpha)$. We can see that $\chi(\omega_X)=\chi(\omega_{\bar{X}})$, thus $\chi(\omega_{\bar{X}})>0$, and so $V^0_{a'}(\omega_{\bar{X}})$ contains a dense open subset of $\Pic^0(A)$. It follows that  the cohomology support locus $$V_{a'}^0(\omega_{\bar{X}}):=\{\alpha\in \Pic^0(A)| H^i(\bar{X},\omega_{\bar{X}}\ot a'^*\alpha)\neq 0\}$$ generates  $\Pic^0(A)$.   By Theorem 2.3 in the paper \cite{chen2001pluricanonical},  we have $\bar{X}$ is of general type.
		\qed\end{proof}
	
	\begin{proof}[Proof of Theorem \ref{thm: bound}]]
		Set $\bar{X}:=X/G$; we perform a $G$-equivalent resolution of the quotient map $\pi\co X\rightarrow\bar{X}$ to obtain the following commutative diagram
		\begin{equation*} 
			\xymatrix{
				Y    \ar[r]^{\bar{\pi}} \ar[d]  &   \bar{Y}  \ar[d]\\
				X \ar[r]^{\pi}  &   \bar{X}        }
		\end{equation*}
		where $\bar{Y}\rightarrow \bar{X}$ is any nonsingular resolution of $\bar{X}$, then $G$ acts on $\bar{Y}\times_{\bar{X}}X$ as $g\cdot(y,x)=(y,g\cdot x)$ for any $y\in\bar{Y},x\in X$; take $Y$ to be the equivariant resolution of singularities on the component of $\bar{Y}\times_{\bar{X}}X$ which dominants $X$ \cite[Theorem 0.1]{AB96}, and we get a generically finite map $\bar{\pi}\co Y\rightarrow \bar{Y}$ of degree $|G|$. By \cite[Lemma 4.3]{Ar08}, we have 
		\begin{equation}\label{eq: Vol Y= G Vol Y+B}
			\Vol(K_Y)=|G|\Vol(K_{\bar{Y}}+B)
		\end{equation} 
		for some effective $\QQ$-divisor $B$ on $\bar{Y}$. Since $X$ has only Gorenstein canonical singularities and $K_X$ is nef, it is a minimal model of $Y$; by Proposition \ref{prop: volume Vol K_X'= Vol K_X}, we have 
		\begin{equation}\label{eq: K^3=Vol Y}
			K_X^3=\mathrm{Vol}(K_Y),
		\end{equation}
		and by  Remark \ref{rem: Vol D leq Vol D+F}, it follows that 
		\begin{equation}\label{eq: K^3Y leq Vol Y+B}
			K^3_{Y_{\min}}=\Vol(K_{\bar{Y}})\leq\Vol(K_{\bar{Y}}+B)
		\end{equation}
		where $Y_{\min}$ is a minimal model of $\bar{Y}$. Combining (\ref{eq: Vol Y= G Vol Y+B}),(\ref{eq: K^3=Vol Y}),(\ref{eq: K^3Y leq Vol Y+B}), $K_X^3\leq 72\chi(\omega_X)$ the Miyaoka-Yau inequality for $X$, and $12\chi(\omega_{\bar{Y}})\leq K_{Y_{\min}}^3$ the Clifford-Severi inequality \cite{Bar15,ZT14} for $Y_{\min}$ (the smooth model $\bar{Y}$ of $Y_{\min}$ is of general type and has maximal Albanese dimension by (2) of Lemma \ref{lem: quotient of X by AutQ}) shows that 
		\begin{equation}\label{eq: compare of volumes}
			12|G|\chi(\omega_{\bar{Y}})\leq|G|K^3_{Y_{\min}}\leq K_X^3\leq72\chi(\omega_X)
		\end{equation}
		Since the map $\pi$ is finite and $X$ has Gorenstein canonical singularities, we can see that $\bar{X}$ has rational singularities (see \cite[Proposition 5.13]{KM98} or \cite[Proposition 1.7]{reid1980canonical}). By Remark \ref{rem: rational singularity} we obtain $\chi(\omega_{\bar{Y}})=\chi(\omega_{\bar{X}})$. From formula (\ref{eq: compare of volumes}) and (\ref{eq: chi X= chi X' = chi X''}) we have $|G|\leq6$. Moreover if $X$ is nonsingular and $K_X$ is ample, replacing the inequality $K_X^3\leq 72\chi(\omega_X)$ by Yau's inequality $K_X^3\leq 64\chi(\omega_X)$ \cite[Remarks: (iii)]{yau1977calabi}, we get $|G|\leq 5$.
		\qed\end{proof}
	
	\begin{remark}\label{rem: bound for threefold isogenous} 
		If $X$ is threefold isogenous to a product of curves, then its invariants satisfies $K_X^3=48\chi(\omega_X)$. Thus in this case we have  $|\mathrm{Aut}_{\mathbb{Q}}(X)|\leq 4$.
	\end{remark}
	\section{Rationally cohomologically rigidity for Albanese general type varieties}\label{sec:4}
	
	We recall that a projective complex variety is said to be of \emph{maximal Albanese dimension} if its Albanese map is generically finite onto its image. According to Catanese \cite{Cat91}, it is said to be of \emph{Albanese general type}, if moreover, its Albanese map is not surjective. We say that a variety $Y$ admits a \emph{higher irrational pencil} if $Y$ admits a surjective morphism with connected fibres onto a nonsingular curve $D$ of genus $g(D)\geq2$. $Y$ is called \emph{Albanese primitive} if it doesn't admit any higher irrational pencil.
	
	Based on the above definitions, we consider the following three classes of projective varieties.
	
	\begin{description}
		\item[Class I] $Y$ is of Albanese general type, and there is some higher irrational pencil $g\co Y\rightarrow D$ whose general fibre $F$ is of Albanese general type.
		\item[Class II] $Y$ is of Albanese general type, and for any higher irrational pencil $g\co Y\rightarrow D$, its general fibre $F$ is of Albanese primitive.
		\item[Class III] $Y$ is of Albanese primitive.	
	\end{description} 
	
	If	 $Y$ belongs to either class I or class II, let $g\co Y\rightarrow D$ be a higher irrational pencil, and let $F$ be its generic fibre.
	Let $\sigma$ be an automorphism of $Y$ such that $g\circ\sigma=g$, we denote the restriction of $\sigma$ on $F$ by $\sigma_F$. The knowledge on the classification of pairs $(F,\sigma_F)$ can help us to understand the classification of pairs $(Y,\sigma)$. For example, $\sigma_F=\mathrm{id}$ implies $\sigma=\mathrm{id}$. Therefore, we have an injective homomorphism
	$$\mathrm{Aut}(Y/D)\rightarrow\mathrm{Aut}(F).$$
	
	The following result shows that if $Y$ admits a fibration over a curve $D$ with $g(D)\geq 1$, then its $\mathrm{Aut}_{\mathbb{Q}}(Y)$ is controlled by the autmorphism group of its generic fibre.
	\begin{lemma}\label{lem: g geq 1 AutQ to Aut F}
		Let $Y$ be a projective variety with $\chi(\omega_Y)>0$, and let $g\co Y\rightarrow D$ be a surjective morphism with connected fibres, where $D$ is a nonsingular curve. Assume that $g(D)\geq 1$. Then we have $\mathrm{Aut}_{\mathbb{Q}}(Y)\subseteq \mathrm{Aut}(Y/D)$. Moreover, if $F$ is a general fibre of $g$, then the induced homomorphism $\mathrm{Aut}_{\mathbb{Q}}(Y)\rightarrow \mathrm{Aut}(F)$ is injective.
	\end{lemma}
	
	\begin{proof} To see $\mathrm{Aut}_{\mathbb{Q}}(Y)\subseteq \mathrm{Aut}(Y/D)$. Take $\sigma\in\Aut_{\QQ}(X)$. Since $\sigma$ acts trivially on $H^2(Y,\mathbb{Q})$, it induces the trivial action on $\mathrm{NS}(Y)\otimes\mathbb{Q}$, where $\mathrm{NS}(Y)$ is the N\'{e}ron-Severi group of $Y$. We conclude that $\sigma(F)$ is numerically equivalent to $F$, where $F$ is any general fibre of $g$, hence $\sigma(F)$ is a fibre of $g$. Let $\bar{\sigma}\co D\rightarrow D$ be the map defined by $\bar{\sigma}(x) = g(\sigma(y))$ for any $y\in g^{-1}(x)$ and $x\in D$. It is sufficient to show that $\bar{\sigma}$ is the identity map. Note that $$g^*\co H^0(D,\Omega^1_{D})\hookrightarrow H^0(Y,\Omega^1_{Y})$$ is an injection, which implies that $\bar{\sigma}$ induces the trivial action on $H^0(D,\Omega^1_{D})$. Therefore, the quotient curve $D/\langle\bar{\sigma}\rangle$ has the same genus as $D$. By Hurwitz formula we obtain
		$$2g(D)-2=o(\bar{\sigma})(2g(D/\langle\bar{\sigma}\rangle)-2 + R)$$
		where $R\geq 0$. It follows that $2g(D)-2\geq o(\bar{\sigma})(2g(D)-2)$. If $g(D)\geq 2$, then we have $ o(\bar{\sigma})=1$, so $\bar{\sigma}=\mathrm{id}$. If $g(D)=1$, then $\bar{\sigma}$ is a translation on $D$. Since $\chi(\omega_Y)>0$, by the proof of Lemma \ref{lem: quotient of X by AutQ}, we have $Y^{\sigma}\neq\varnothing$, hence $\bar{\sigma}$ has fixed pints and $\bar{\sigma}=\mathrm{id}$.
		
		Let $F$ be a general fibre of $g$. Composing the inclusion $\mathrm{Aut}_{\mathbb{Q}}(Y)\hookrightarrow \mathrm{Aut}(Y/D)$ with the injective homomorphism $\mathrm{Aut}(Y/D)\rightarrow \mathrm{Aut}(F)$, we get an injective homomorphism $\mathrm{Aut}(Y/D)\rightarrow\mathrm{Aut}(F)$.
		\qed
	\end{proof}
	
	Recall that for a nonsingular projective surface $S$ which admits a fibration $f\co S\to B$ of genus $g\geq 2$, if there is a nontrivial automorphism $\sigma$ in $\Aut(S/B)$ which induces a trivial action on $H^0(S,\omega_S)$, then $g(B)\leq 1$. This was proved by Cai \cite[Lemma 2.1]{Cai12}. The following result is a generalization of this to higher dimension. It shows that if $Y$ admits a higher irrational pencil, then its $\mathrm{Aut}_{\mathbb{Q}}(Y)$ is controlled by the subgroup of automorphisms of its general fibre $F$ acting trivially on $H^0(F,\omega_F)$.
	
	\begin{lemma}\label{lem: g geq 2 AutQ to Aut F} 
		Let $g\co Y\rightarrow D$ be a higher irrational pencil of a projective variety $Y$, and $F$ be a general fibre of $g$. Let $\sigma\in\mathrm{Aut}(Y/D)$.
		
		If $\sigma$ induces a trivial action on $H^0(Y,\omega_Y)$, then $\sigma_F$ induces identity on $H^0(F,\omega_F)$. Moreover, if in addition $\chi(\omega_Y)>0$, then $\mathrm{Im}(\mathrm{Aut}_{\mathbb{Q}}(Y)\rightarrow \mathrm{Aut}(F))$ consists of automorphisms of $F$ acting trivially on $H^0(F,\omega_F)$.
	\end{lemma}
	\begin{proof} Consider the induced action of $\sigma$ on $g_*\omega_Y$, we can decomposite it as $g_*\omega_Y=\mathcal{E}\oplus\mathcal{F}$ with eignsubseaf $\mathcal{E}$ of eignvalue $=1$ and direct sum of eignsubsheaves $\mathcal{F}$ of eignvalues $\neq 1$. Let $\mathcal{E}'\subset g_*\omega_Y$ be the subsheaf generated by global sections of $g_*\omega_Y$. Since $\sigma$ acts trivially on $H^0(Y,\omega_Y)$, we have $\mathcal{E}'\subseteq \mathcal{E}$. Therefore, $h^0(D,\mathcal{E})=h^0(D,g_*\omega_Y)$ and hence $h^0(D,\mathcal{F})=0$. Applying the Riemann-Roch formula to $\mathcal{F}$ we obtain
		$$\deg(\mathcal{F})+r(1-g(D))=-h^1(D,\mathcal{F})\leq 0.$$
		The semi-positivity of $g_*\omega_Y\otimes\omega^{-1}_D$ imples that $$\deg(\mathcal{F})-2r(g(D)-1)=\deg(\mathcal{F}\otimes \omega^{-1}_D)\geq 0$$ By the assumption that $g\co Y\rightarrow D$ is a higher irrational pencil, we get $g(D)\geq 2$. Combining the two above inequalities, we have $r=0$ and hence $\mathcal{F}=0$. Note that the natrual map $g_*\omega_Y\otimes\mathbb{C}(p)\rightarrow H^0(F,\omega_F)$ is an isomorphism, where $F=g^{-1}(p)$ for a general point $p\in D$, it follows that $\sigma_F$ induces a trivial action on $H^0(F,\omega_F)$. 
		
		By the proof of Lemma \ref{lem: g geq 1 AutQ to Aut F}, $g(D)\geq 2$ implies that $\mathrm{Aut}_{\mathbb{Q}}(Y)\subseteq\mathrm{Aut}(Y/D)$. Hence for all $\sigma\in\mathrm{Aut}_{\mathbb{Q}}(Y)$, $\sigma_F$ induces the trivial action on $H^0(F,\omega_F)$.
		\qed\end{proof}
	
	To apply induction, we now consider the following situation:
	\begin{description}
		\item[Situation (*)]Let $Y$ be a projective variety of $\dim Y\geq 3$. Suppose there are a seqence higher irrational pencils $g_j\co F_j\rightarrow D_j$ such that each $F_{j+1}$ is a general fibre of $g_j$ for $0\leq j\leq \dim Y-2$, set $F_0=Y$. We define
		$$\mathcal{W}_1:=\mathrm{Im}(\mathrm{Aut}(F_0/D_0)\rightarrow\mathrm{Aut}(F_1))\cap\mathrm{Aut}(F_1/D_1)$$
		and for $2\leq j\leq \dim Y-2$
		$$\mathcal{W}_k:=\mathrm{Im}(\mathcal{W}_{k-1}\rightarrow\mathrm{Aut}(F_k))\cap\mathrm{Aut}(F_k/D_k).$$
		For any $\sigma\in \mathcal{W}_k$, it uniquely determines a sequence of automorphisms $\tilde{\sigma}_j\in\mathrm{Aut}(F_j/D_j)$ for $0\leq j\leq k-1$ such that $\tilde{\sigma}_{j+1}=\tilde{\sigma}_j|_{F_{j+1}}$, where $\tilde{\sigma}_k=\sigma$. We call $\sigma$ the \emph{successive restriction} of $\tilde{\sigma}_0$.
	\end{description}
	
	\begin{corollary}\label{cor: Situation *}
		Let $Y$ be a projective variety as in Situation (*). Then there is no nontrivial automorphism $\sigma\in \mathcal{W}_{\dim Y-2}$ such that $\tilde{\sigma}_0$ induces trivial action on $H^0(Y,\omega_Y)$.
	\end{corollary}
	
	\begin{proof} Set $d:=\dim Y$. By the definition of $\mathcal{W}_{\dim Y-2}$, there is a sequence of automorphisms $\tilde{\sigma}_j\in\mathrm{Aut}(F_j/D_j)$ for $0\leq j\leq d-3$ such that $\tilde{\sigma}_{j+1}=\tilde{\sigma}_j|_{F_{j+1}}$, where $\tilde{\sigma}_{d-2}=\sigma$. Applying Lemma \ref{lem: g geq 2 AutQ to Aut F} repeatly to fibration $g_j$, we see that $\tilde{\sigma}_{j+1}$ induces the trivial action on $H^0(F_{j+1},\omega_{F_{j+1}})$ for $0\leq j\leq d-3$. In particular, $\tilde{\sigma}_{d-2}$ induces the trivial action on $H^0(F_{d-2},\omega_{F_{d-2}})$ and it belongs to $\mathrm{Aut}(F_{d-2}/D_{d-2})$. Note that $F_{d-2}$ is a nonsingular surface and $g(D_{d-2})\geq 2$, from Lemma 2.1 of \cite{Cai12}, we conclude that $\sigma=\tilde{\sigma}_{d-2}=\mathrm{id}$. 
		\qed\end{proof}
	
	\begin{example}\label{ex: d fold isogenous with a sequence of higher irrational pencil}
		Let $X$ be a $d$-fold isogenous to an unmixed product of curves, and let $$(C_1\times\dots\times C_d)/G$$ be its minimal realization. Suppose $X$ belongs to class I. Since $X$ is of Albanese general type and $q(X)=\sum_{i=1}^{d}g(C_i)$, we have $g(\bar{C}_j)\geq 1$ for all $1\leq j\leq d$ and $g(\bar{C}_i)\geq 2$ for some $1\leq i\leq d$. Without loss of generality, suppose $g(\bar{C}_d)\geq 2$, then $f_d\co X\rightarrow\bar{C}_d$ is a higher irrational pencil. Let $U_d$ be the  complement of the branch points of the quotient map $C_d\rightarrow\bar{C}_d$, then for all $x\in U_d$, $F_x$ is of Albanese general type.
		
		Since $F_x$ is isomorphic to $(C_1\times\dots \times C_{d-1})/K_d$ which is a $d-1$-fold isogenous to an unmxied product of curves.  For the same reason, $F_x$ admits a higher irrational pencil. Without loss of generality, suppose $g(C_{d-1}/K_d)\geq 2$, hence $g\co F_x\rightarrow C_{d-1}/K_d$ is the corresponding higher irrational pencil. Repeat this process, we get a sequence of higher irrational pencils as following:
		\begin{equation}
			\xymatrix{
				F_{d-2}\ar[r] \ar[d]^{g_{d-2}} & \dots \ar[r] & F_2 \ar[r] \ar[d]^{g_2} & F_1 \ar[r] \ar[d]^{g_1} & X \ar[d]^{f_d} \\
				D_{d-2}  & \dots  & D_2  & D_1  & \bar{C}_d
			}
		\end{equation}
		where $F_1\cong (C_1\times \dots \times C_{d-1})/K_d$ and $D_1=C_{d-1}/K_d$. Each $F_{j+1}$ is a general fibre of the higher irrational pencil $g_i$. Therefore, $F_2$ is isomorphic to $$(C_1\times\dots\times C_{d-2})/K_d\cap K_{d-1}=C_1\times\dots\times C_{d-2}$$ provided that $K_d\cap K_{d-1}=\{1\}$, and we have $F_j=C_1\times\dots\times C_{d-j}$ and $D_j=C_{d-j}$ for $2\leq j\leq d-2$. Therefore $X$ is in Situation (*).
	\end{example}

	\begin{theorem}\label{thm: d-folds isogenous AutQ}
		Let $X$ be $d$-fold isogenous to an unmixed product of curves with $d\geq 3$, and $(C_1\times \dots\times C_d)/G$ be its minimal realization. Suppose $g(C_i/G)\geq 1$ for all $1\leq i\leq d$.  If there is some $1\leq i\leq d$ with $g(C_i/G)\geq 2$ such that $g(C_j/K_i)\geq 2$ for some $j\neq i$, then $\mathrm{Aut}_{\mathbb{Q}}(X)$ is trivial. 
	\end{theorem} 
	
	\begin{proof} By Example \ref{ex: d fold isogenous with a sequence of higher irrational pencil}, $X$ is in Situation (*), we get a sequence of higher irrational pencils $g_j\co F_j\rightarrow D_j$ for $1\leq j\leq d-2$. Set $F_0=X,D_0=\bar{C}_d$ and $g_0=f_d$. We may assume \begin{align*}
			&F_1\cong (C_1\times\dots \times C_{d-1})/K_d,D_1=C_{d-1}/K_d \text{ and }\\
			&F_j\cong C_1\times\dots \times C_{d-j}, D_j=C_{d-j}\text{ for } j\geq 2,
		\end{align*} each $F_{j+1}$ is a general fibre of $g_j$. 
		
		Let $\tau\in\mathrm{Aut}_{\mathbb{Q}}(X)$, and let $k$ be the maximal integer such that $\sigma\in\mathcal{W}_k$ which is the successive restriction of $\tau$. If $k=d-2$. Since $\tau$ acts trivially on $H^0(X,\omega_X)$, by Corollary \ref{cor: Situation *}, $\sigma=\mathrm{id}$ and hence $\tau=\mathrm{id}$.
		
		Now assume that  $\tau\neq \mathrm{id}$ and $k<d-2$. Since $\tau$ induces the trivial action on $H^0(F_0,\omega_{F_0})$, by Lemma \ref{lem: g geq 2 AutQ to Aut F}, we know that  $\tilde{\sigma}_j$ induces identity on $H^0(F_j,\omega_{F_j})$ and $\tilde{\sigma}_j\in\mathrm{Aut}(F_j/D_j)$ for $1\leq j\leq k$. In particular, we have $\sigma'=\tilde{\sigma}_k|_{F_{k+1}}$ induces the trivial action on $H^0(F_{k+1},\omega_{F_{k+1}})$ and $\sigma'\notin\mathrm{Aut}(F_{k+1}/D_{k+1})$.
		
		\noindent\textbf{Step 1.} \emph{$\sigma'$ induces an automorphism of $D_{k+1}$}. 
		
		Let $U_{k+1}'$ and $U_{k+1}$ be Zariski open subsets of $D_{k+1}$ and $\bar{C}_{d-k-1}=C_{d-k-1}/G$, respectively, such that the corresponding map $\rho\co U_{k+1}'\rightarrow U_{k+1}$ is unramified. Set $W_y:=g_{k+1}^{-1}(y)$ where $y\in U_{k+1}'$ and $x=\rho(y)$. Let $f_{d-k-1}\co X\rightarrow \bar{C}_{d-k-1}$ be the fibration induced by the natrual projection, then there is a commutative diagram
		\begin{equation}\label{eq: commutative diagram F k+1 to X}
			\xymatrix{
				F_{k+1} \ar[r]^{\eta} \ar[d]^{g_{k+1}} & X\ar[d]^{f_{d-k-1}}\\
				D_{k+1}\ar[r]^{\rho} & \bar{C}_{d-k-1}
			}
		\end{equation}
		where $\eta$ is the natrual embedding. By Lemma \ref{lem: g geq 1 AutQ to Aut F}, $f_{d-k-1}\circ\tau=f_{d-k-1}$ and hence 
		$$f_{d-k-1}(\eta(\sigma'(W_y)))=f_{d-k-1}(\tau(\eta(W_y)))=f_{d-k-1}(\eta(W_y))=\rho(y)=x$$
		It follows that $\sigma'(W_y)$ is a fibre of $g_{k+1}$ for any $y\in U_{k+1}'$, then $\sigma'$ induces an automorphism $\bar{\sigma}$ of $D_{k+1}$. 
		
		\noindent\textbf{Step 2.} \emph{$\rho$ factors through the quotient map $D_{k+1}\rightarrow D_{k+1}/\langle\bar{\sigma}\rangle$}.
		
		By assumption $\bar{\sigma}\neq \mathrm{id}$, otherwise $g_{k+1}\circ\sigma'=g_{k+1}$. Consider the induced action of $\bar{\sigma}$ on $H^0(D_{k+1},\omega_{D_{k+1}})$, we have a decomposition $$H^0(D_{k+1},\omega_{D_{k+1}})=V\oplus V',$$ where $V$ is the subspace with eignvalue $=1$, $V'$ is the direct sum of subspaces with eignvalues $\neq 1$. By the definition of $\bar{{\sigma}}$ we know that: for any $x\in D_{k+1}$, $\sigma'(W_x)=W_{\bar{\sigma}(x)}$ where $W_x=g_{k+1}^{-1}(x)$. Since $f_{d-k-1}\circ\tau = f_{d-k-1}$ we have
		\begin{equation}\label{eq: f eta sig W = f eta W}
			f_{d-k-1}(\eta(W_{\bar{\sigma}(x)}))=f_{d-k-1}(\eta(\sigma'(W_x)))=f_{d-k-1}(\tau(\eta(W_x))=f_{d-k-1}(\eta(W_x)).
		\end{equation}
		The commutative diagram \ref{eq: commutative diagram F k+1 to X} shows that \begin{equation*}
			\text{ for any } y\in W_x\text{ we have }f_{d-k-1}(\eta(y))=\rho(g_{k+1}(y))=\rho(x).
		\end{equation*} From the formula (\ref{eq: f eta sig W = f eta W}) we know that 
		$$f_{d-k-1}(\eta(y))=f_{d-k-1}(\eta(\sigma'(y)))=\rho(\bar{\sigma}(x)),$$
		hence that $\rho\circ\bar{{\sigma}}=\rho$. Set $\bar{D}_{k+1}=D_{k+1}/\langle\bar{\sigma}\rangle$, hence the morphism $\rho$ factors through the quotient map $D_{k+1}\rightarrow \bar{D}_{k+1}$.
		
		\noindent\textbf{Step 3.}
		
		By step 2, $g(\bar{D}_{k+1})\geq g(\bar{C}_{k+1})\geq 1$ and $\dim V \geq 1$. Since $\bar{{\sigma}}\neq \mathrm{id}$, $\dim V'\geq 1$. By step 1, $\sigma'$ induces the trivial action on 
		$H^0(F_{k+1},\omega_{k+1})$. Note that $F_{k+1}\cong C_1\times \dots\times C_{d-k-1}$ for $ k>1$ and $F_1\cong (C_1\times\dots\times C_{d-1})/K_d$. Let $h_j\co F\rightarrow D_j$ be the fibration induced by the natrual projection, then we have an injection $$h_1^*H^0(D_1,\omega_{D_1})\wedge\dots\wedge g_{k+1}^*H^0(D_{k+1},\omega_{D_{k+1}})\rightarrow H^0(F_{k+1},\omega_{F_{k+1}})$$
		Therefore $\sigma'$ acts trivially on $h_1^*H^0(D_1,\omega_{D_1})\wedge\dots\wedge g_{k+1}^*H^0(D_{k+1},\omega_{D_{k+1}})$. Since the induced action of $\sigma'$ on $g_{k+1}^*H^0(D_{k+1}\omega_{D_{k+1}})$ can be identified as the induced action of $\bar{{\sigma}}$ on $H^0(D_{k+1},\omega_{D_{k+1}})$, the induced action of $\sigma'$ on $g_{k+1}^*H^0(D_{k+1},\omega_{D_{k+1}})$ has two different eign-subsapces, which contradicts to the fact that $\sigma'$ induces identity on $$h_1^*H^0(D_1,\omega_{D_1})\wedge\dots\wedge g_{k+1}^*H^0(D_{k+1},\omega_{D_{k+1}}).$$ Hence $\tau=\mathrm{id}$.\qed
	\end{proof}
	
	\section{ $\mathrm{Aut}_{\mathbb{Q}}(X)$ for threefolds isogenous to an unmixed product of curves}\label{sec:5}
	
	In this section, we focus on the threefolds isogenouse to an unmixed product of curves with maximal Albanese dimension. Our main main result is the following.
	
	\begin{theorem}\label{thm: threefold isogenous AutQ}
		Let $X$ be a threefold isogenous to an unmixed product of curves, and let $(C_1\times C_2\times C_3)/G$ be its minimal realization. Suppose $g(C_i/G)\geq 1$ for all $1\leq i\leq 3$. The kernel of group homomorphism $\psi_i\co G\rightarrow \mathrm{Aut}(C_i)$ will be denoted by $K_i$. Then the following cases occur
		\begin{enumerate} 
			\item if there is a pair $(i,j)$ with $j\neq i$ such that $g(C_i/G)\geq 2$ and $g(C_j/K_i)\geq 2$, then $\mathrm{Aut}_{\mathbb{Q}}(X)$ is trivial;
			\item if for any $1\leq i\leq 3$ with $g(C_i/G)\geq 2$, $g(C_j/K_i)=1$ for all $j\neq i$, then $\mathrm{Aut}_{\mathbb{Q}}(X)\cong(\mathbb{Z}_2)^k$ where $k=0,1$;
			\item if for all $1\leq i\leq 3$, $g(C_i/G)=1$, and suppose that the group $G$ is an abelian group and $K_i$ is a cyclic group for all $1\leq i\leq 3$, then  $\mathrm{Aut}_{\mathbb{Q}}(X)\cong (\mathbb{Z}_2)^k$ where $k=0,1,2$.
	\end{enumerate}  \end{theorem}
	
	\subsection{The case of Albanese general type}\label{sec:5.1}
	Let $X$ be a threefold isogenous to an unmixed product of curves, and let $(C_1\times C_2\times C_3)/G$ be its minimal realization. Suppose $g(C_i/G)\geq 1$ for all $1\leq i\leq 3$ and there is some $1\leq i\leq 3$ such that $g(C_i/G)\geq 2$. Let $F$ be a general fibre of the fibration $f_i\co X\rightarrow \bar{C}_i$ induced by the natrual projection. Fix an element $\sigma\in\mathrm{Aut}_{\mathbb{Q}}(X)$. Then the following properties are satisfied:
	\begin{enumerate}
		\item $X$ is of Albanese general type and $f_i$ is a higher irrational pencil.
		\item $F$ is a surface isogenous to an unmixed product of curves with maximal Albanese dimension.
		\item $f_i\circ\sigma=f_i$ ( Lemma \ref{lem: g geq 1 AutQ to Aut F}).
		\item $\sigma$ induces trivial action on $H^0(F,\omega_F)$ (Lemma \ref{lem: g geq 2 AutQ to Aut F}).
	\end{enumerate}
	
	According to the classification of projective varieties of general type of maximal Albanese dimension in Section \ref{sec:4}. We can divide $X$ into three classes.
	
	\begin{description}
		\item[Class I] There is a pair $(i,j)$ with $j\neq i$ such that $g(C_i/G)\geq 2$ and $g(C_j/K_i)\geq 2$.
		\item[Class II]  For any $1\leq i\leq 3$ with $g(C_i/G)\geq 2$, $g(C_j/K_i)=1$ for all $j\neq i$
		\item[Class III] For all $1\leq i\leq 3$, $g(C_i/G)=1$.
	\end{description} 
	
	\begin{corollary}\label{cor: class I} 	
		If $X$ belongs to class I, then $\mathrm{Aut}_{\mathbb{Q}}(X)$ is trivial.
	\end{corollary}
	
	\begin{proof}
		By Theorem \ref{thm: d-folds isogenous AutQ}, the conclusion holds.\qed
	\end{proof}
	
	\begin{corollary}\label{cor: class II}
		If $X$ belongs to class II, then $|\mathrm{Aut}_{\mathbb{Q}}(X)|\leq 2$.
	\end{corollary}
	
	\begin{proof}
		By the definition of class II, we can assume $g(C_3/G)\geq 2$, then the fibration $f\co X\rightarrow \bar{C}_3=C_3/G$ induced by the natrual projection is a higher irrational pencil, whose general fibre $F\cong(C_1\times C_2)/K_3$ is Albanese primitive. It follows that $g(C_i/K_3)=1$ for $i=1,2$. By Lemma \ref{lem: g geq 1 AutQ to Aut F}, for any $\sigma\in\mathrm{Aut}_{\mathbb{Q}}(X)$, $f\circ\sigma=f$.
		
		\noindent\textbf{Step 1.}\emph{ We claim that $\sigma_F:=\sigma|_F$ induces trivial action on $H^0(F,\Omega^1_F)$.}
		
		Note that 
		\begin{align*}
			&H^0(F,\Omega^1_F)=H^0(C_1,\Omega^1_{C_1})^{K_3}\oplus H^0(C_2,\Omega^1_{C_2})^{K_3}\\
			&H^0(X,\Omega^1_X)=H^0(C_1,\Omega^1_{C_1})^G\oplus H^0(C_2,\Omega^1_{C_2})^G\oplus H^0(C_3,\Omega^1_{C_3})^G
		\end{align*}
		From $1=g(C_i/K_3)\geq g(C_i/G)\geq 1$ we know $g(C_i/K_3)=g(C_i/G)=1$, hence $$\dim H^0(C_i,\Omega^1_{C_i})^{K_3}=\dim  H^0(C_i,\Omega^1_{C_i})^G$$ for $i=1,2$.  Let $j\co F\rightarrow X$ be the natrual embedding. Based on the equalities of dimensions mentioned above, we know that $$j^*\co H^0(X,\Omega^1_X)\rightarrow H^0(F,\Omega^1_F)$$ is a surjective homomorphism. Consider the induced actions of $\sigma$ and $\sigma_F$ on $H^0(X,\Omega^1_X)$ and $H^0(F,\Omega^1_F)$, respectively, there is a commutative diagram
		$$\xymatrix{
			H^0(X,\Omega^1_X) \ar[r]^{j^*} \ar[d]^{\sigma^*} & H^0(F,\Omega^1_F) \ar[d]^{\sigma_F^*}\\
			H^0(X,\Omega^1_X) \ar[r]^{j^*}  & H^0(F,\Omega^1_F)	
		}$$
		For any $v\in H^0(F,\Omega^1_F)$ there exists $u\in H^0(X,\Omega^1_X)$ such that $v=j^*u$, hence $$\sigma_F^*v=\sigma_F^*j^*u=j^*\sigma^*u=j^*u=v.$$ So $\sigma_F^*$ acts trivially on $H^0(F,\Omega^1_F)$. According to Lemma \ref{lem: g geq 2 AutQ to Aut F}, $\sigma_F^*$ induces trivial action on $H^0(F,\omega_F)$. Then we have an injective group homomorphism 
		$$\mathrm{Aut}_{\mathbb{Q}}(X)\hookrightarrow\mathrm{Aut}_{d}(F),\quad\sigma\mapsto\sigma_F,$$
		where $\mathrm{Aut}_{d}(F)$ is the subgroup of automorphisms of $F$ acting trivially on $H^0(F,\Omega^j_F)$ for all $j=1,2$.
		
		\noindent\textbf{Step 2.} We claim that $|\mathrm{Aut}_{d}(F)|\leq 2$.
		
		Note that $F$ is a minimal surface of general type and of maximal Albanese dimension. Using the same method of Lemma \ref{lem: quotient of X by AutQ}, we can show that the Albanese map $a_F\co F\rightarrow A_F$ factors through the the quotient map $F\rightarrow F'=:F/\mathrm{Aut}_{d}(F)$ and $\chi(\omega_F)=\chi(\omega_F')>0$. Let $S$ be the minimal smooth model of $F'$, by G. Xiao's result $|\mathrm{Aut}_{d}(F)|K_S^2\leq K_F^2$ \cite{Xiao}. Combining this with the Severi inequality $4\chi(\omega_S)\leq K_S^2$, the Bogomolov-Miyaoka-Yau inequality $K_F^2\leq 9\chi(\omega_F)$ and $\chi(\omega_F)=\chi(\omega_{F'})=\chi(\omega_S)$, we conclude that $|\mathrm{Aut}_{d}(F)|\leq 2$, and hence $|\mathrm{Aut}_{\mathbb{Q}}(X)|\leq 2$.\qed
	\end{proof}
	
	\subsection{Automorphisms of $X$ descended from $\Aut(\curProd)$}\label{sec:5.2} Throughout the rest of this section, we assume that $g(C_i/G)=1$ for all $1\leq i\leq 3$ and $G$ is an abelian group. Set $\mathrm{Lf}(X)=(G\times G\times G)/K\Delta_G$ where $K:=K_1\times K_2\times K_3$, $K_i=\mathrm{Ker}(G\rightarrow\mathrm{Aut}(C_i))$ and $\Delta_G$ the diagonal subgroup of $G\times G\times G$.
	\begin{lemma}\label{lem: LfX to AutX}
		With the above notations, we have an injective group homomorphism $$\mathrm{Lf}(X)\hookrightarrow \Aut(X).$$
	\end{lemma}
	
	\begin{proof}Recall that homomorphism $\psi_i\co G\rightarrow\Aut(C_i)$ is given by the $G$-action on $C_i$ for each $1\leq i\leq 3$, then we get a homomorphism $$\psi:=\psi_1\times\psi_2\times\psi_3\co G\times G\times G\rightarrow \Aut(\curProd).$$ Since $K_i=\Ker(\psi_i)$ for $1\leq i\leq 3$, we can see that $K=\Ker(\psi)$. As $X$ is the quotient of $\curProd$ under the action of $\psi(\Delta_G)$, it follows that
		\begin{equation}\label{eq: isomorphism of AutX to N/D}
			\mathrm{Aut}(X)\cong N(\psi(\Delta_G))/\psi(\Delta_G),
		\end{equation} 
		where $N(\psi(\Delta_G))$ the normalizer of $\psi(\Delta_G)$ in $\Aut(\curProd)$. Since $G$ is abelian and $\psi$ is a group homomorphism,
		\begin{equation*}
			\quad\psi(h g)=\psi(g)\psi(h), \text{ for any } g\in\Delta_G\text{ and }h\in G\times G\times G, 
		\end{equation*} hence $\psi(G\times G\times G)\subseteq N(\psi(\Delta_G))$.  Composing $\psi$ with the qoutient map $N(\psi(\Delta_G))\rightarrow N(\psi(\Delta_G))/\psi(\Delta_G)$ and the isomorphism (\ref{eq: isomorphism of AutX to N/D}), we have a group homomorphism $$G\times G\times G\rightarrow\Aut(X)$$ with kernel $K\Delta_G$, which induces an injective homomorphism $$(G\times G\times G)/K\Delta_G\hookrightarrow\mathrm{Aut}(X).$$
		\qed\end{proof}
	
	Note that we have an injective homomorphism $$j\co  (G\times G\times G)/K\hookrightarrow \Aut(\curProd)$$ and that $((G\times G\times G)/K)/\Delta_G=\mathrm{Lf}(X)$. For any $g\in\mathrm{Lf}(X)$, there exists some $\tilde{g}\in (G\times G\times G)/K$ such that $\tilde{g}\Delta_G=g$, we call the image $j(g)$ of $\tilde{g}$ in $\mathrm{Aut}(\curProd)$ a lifting of $g$. The following proposition shows that numerically trivial automorphism of $X$ can be lifted to $\Aut(\curProd)$.
	
	\begin{proposition}\label{prop: AutQ to LfX} 
		Let $\sigma$ be an automorphism of $X$ satisfies the property: for all $1\leq i\leq 3$,  $f_i\circ\sigma=f_i$, where $f_i\co X\rightarrow C_i/G$ is a fibration induced by the natural projtection, then there is an injective group homomorphism $\langle\sigma\rangle\hookrightarrow\mathrm{Lf}(X)$.
		
		In particular, we have an injective homomorphism $$\mathrm{Aut}_{\mathbb{Q}}(X)\hookrightarrow\mathrm{Lf}(X).$$ \end{proposition}
	
	\begin{proof} Set $\bar{C}_i=C_i/G$. Since $(\curProd)/G$ is the minimal realization of $X$, we have $K\cap\Delta_G=\{1\}$, and we can identify $\Delta_G$ as its image in the quotient group $(G\times G\times G)/K$ under the quotient map. Therefore, we have an injective homomorphism $$\Delta_G\hookrightarrow G/K_1\times G/K_2 \times G/K_3.$$ Consider the action of $G/K_1\times G/K_2 \times G/K_3$ on $C_1\times C_2\times C_3$, there is a commutative diagram:
		\begin{equation*}
			\xymatrix{
				C_1\times C_2\times C_3\ar[rr]^{\pi} \ar[rd]^{\varphi} & & X \ar[ld]^{\varphi'}\\
				& \bar{C}_1\times\bar{C}_2\times\bar{C}_3 &
			}
		\end{equation*}
		where $\bar{C}_1\times\bar{C}_2\times\bar{C}_3=(C_1\times C_2\times C_3)/G/K_1\times G/K_2 \times G/K_3$. Hence \begin{align*}
			&\mathrm{Gal}(\varphi)\cong G/K_1\times G/K_2 \times G/K_3,\\
			&\mathrm{Gal}(\pi)\cong\Delta_G.
		\end{align*}It follows that $$\mathrm{Gal}(\varphi')=\mathrm{Gal}(\varphi)/\mathrm{Gal}(\pi)\cong (G/K_1\times G/K_2 \times G/K_3)/\Delta_G\cong \mathrm{Lf}(X).$$
		Set $H:=\langle \sigma\rangle $ and $Y:=X/H$. By assumption, $f_i\circ\sigma=f_i$ for all $1\leq i\leq 3$. Note that $\varphi'=f_1\times f_2\times f_3$, it follows that $\varphi'$ factors through the quotient map $\varphi''\co X\rightarrow Y$, \begin{equation*}
			\xymatrix{
				X \ar[rr]^{\varphi''} \ar[rd]^{\varphi'} & & Y \ar[ld]\\
				& \bar{C}_1\times\bar{C}_2\times\bar{C}_3 &
			}
		\end{equation*}
		Theorefore $H\cong\mathrm{Gal}(\varphi'')\triangleleft\mathrm{Gal}(\varphi')\cong\mathrm{Lf}(X)$, and thus there is an injective homomorphism $H\hookrightarrow\mathrm{Lf}(X)$.
		
		By lemma \ref{lem: g geq 1 AutQ to Aut F}, take any $\sigma\in\mathrm{Aut}_{\mathbb{Q}}(X)$, $f_i\circ\sigma=f_i$ for all $1\leq i\leq 3$. According to the argument above, there is an injective homomorphism $\mathrm{Aut}_{\mathbb{Q}}(X)\hookrightarrow\mathrm{Lf}(X)$.\qed
	\end{proof}

	\begin{remark}\label{rem: example of AutQ can't lifts} 
		For irregular surfaces which is not of maximal Albanese dimension, Cai and Liu find one $S$ surface isogenous to a product of curves with $q(S)=1$ and $\mathrm{Aut}_{\mathbb{Q}}(S)\cong\ZZ_4$ \cite[Example 4.6]{cai2018automorphisms}, a generator of this group can't lift to an automorphism of the product of curves associated to the minimal realization of  $S$.   
	\end{remark}
	
	\subsubsection{Representations of $\mathrm{Lf}(X)$}\label{sec:5.3} Form now to the end of this section, we fix an algebaric data $$\mathcal{A}=(G,K_1,K_2,K_3,\cV_1,\cV_2,\cV_3)$$ for some threefold $X$ isogenous to an unmixed product of curves. 
	\begin{definition}\label{def: admissible characters}
		A linear character $\chi_1\times\chi_2\times\chi_3$ of group $G\times G\times G$ is called \emph{admissible character} for $\mathcal{A}$ if it satisfies the following conditions:
		\begin{enumerate}[(1)]
			\item $K_i\subseteq\Ker(\chi_i)$ for all $1\leq i\leq 3$;
			\item if $\chi_i\neq 1_G$, then $\chi_i(\sigma)\neq 1$ for some $\sigma\in \Sigma_i$;
			\item $\chi_1\chi_2\chi_3=1_G$.
		\end{enumerate}The number of characters $\chi_i$ such that $\chi_i\neq 1_G$ of an admissible character $\chi_1\times\chi_2\times\chi_3$ is called the \emph{weight}. Denote the set of all admissible characters of weight $3$ (resp. weight $2$)  by $A_3$ (resp. $A_2$), and set $A=A_3\cup A_2$. 
	\end{definition}
	Let $\chi_1\times\chi_2\times\chi_3$ be an admissible characater for $\mathcal{A}$. Recall that $\mathrm{Lf}(X)=(G\times G\times G)/K\Delta_G$. Since the conditions (1) and (3) in Definition \ref{def: admissible characters} implies that $K\Delta_G\subseteq\Ker(\chi_1\times\chi_2\times\chi_3)$, the admissible character $\chi_1\times\chi_2\times\chi_3$ can be regard as a linear character of the group $\mathrm{Lf}(X)$. We next consider the actions of $\mathrm{Lf}(X)$ on $H^i(X,\CC)$ for $1\leq i\leq 3$ given by
	$$\rho_i\co \mathrm{Lf}(X)\rightarrow \GL(H^i(X,\CC)),\quad g\mapsto [\omega\mapsto(g^{-1})^*\omega].$$ 
	
	\begin{lemma}\label{lem: representation of LfX}
		Under the above notations. We have the following properties:
		\begin{enumerate}
			\item The representation $\rho_1$ is trivial;
			\item and for any $\psi\in\mathrm{Lf}(X)^*$ and $i=2,3$, the character space $H^i(X,\mathbb{C})^{\psi}$ under the representation $\rho_i$ is non-zero  if and only if $\psi$ is an admissible character. 
		\end{enumerate}  Moreover, there is a filtration of subgroups of $\mathrm{Lf}(X)$
		$$\Ker(\rho_3)\subset\Ker(\rho_2)\subset\Ker(\rho_1)=\mathrm{Lf}(X).$$
		Therefore, $\mathrm{Aut}_{\mathbb{Q}}(X)=\Ker(\rho_3)=(\cap_{\psi\in A}\Ker(\psi))/K\Delta_G.$
	\end{lemma}
	\begin{proof} Restricting the $G$-action of $\curProd$ on $C_i$ we get representations $$\varphi_i\co G\rightarrow\GL(H^1(C_i,\CC)).$$ By Theoerm \ref{thm: Chevalley-Weil formula} we have $$H^1(C_i,\CC)=H^1(C_i,\CC)^{1_G}\oplus(\oplus_{\chi\in I_i}H^1(C_i,\CC)^{\chi})$$ where $I_i=\{\chi\in G^*|\chi(\sigma)\neq1 \text{ for some }\sigma\in\Sigma_i \text{ and }K_i\subset\Ker(\chi)\}$. 
		
		By the K\"{u}nneth theorem of the cohomology of product spaces, we can see that 
		\begin{align*}H^1(C_1\times C_2\times C_3,\mathbb{C})&= \bigoplus_{\substack{1\leq i<j\leq3,\\ k\neq i,j}}H^1(C_i,\mathbb{C})\otimes H^0(C_j,\mathbb{C})\otimes H^0(C_k,\mathbb{C})\\
			H^2(C_1\times C_2\times C_3,\mathbb{C})&=W_2\bigoplus\bigg( \bigoplus_{\substack{1\leq i<j\leq3,\\ k\neq i,j}}H^1(C_i,\mathbb{C})\otimes H^1(C_j,\mathbb{C})\otimes H^0(C_k,\mathbb{C})\bigg)\\
			H^3(C_1\times C_2\times C_3,\mathbb{C})&=W_3\bigoplus\bigg(H^1(C_1,\mathbb{C})\otimes H^1(C_2,\mathbb{C})\otimes H^1(C_3,\mathbb{C})\bigg)
		\end{align*}
		where \begin{align*}
			&W_2=\bigoplus_{\substack{1\leq i<j\leq3,\\ k\neq i,j}}H^2(C_i,\mathbb{C})\otimes H^0(C_j,\mathbb{C})\otimes H^0(C_k,\mathbb{C}),\\
			&W_3=\bigoplus_{1\leq i,j\leq3}H^2(C_i,\mathbb{C})\otimes H^1(C_j,\mathbb{C})\otimes H^0(C_k,\mathbb{C}).\end{align*}
		Since $H^i(X,\CC)\cong H^i(\curProd,\CC)^{G}$, we obtain following decompositions:
		\begin{align}
			H^1(X,\mathbb{C})&=\bigoplus_{\substack{1\leq i<j\leq3,\\ k\neq i,j}}H^1(C_i,\mathbb{C})^{1_G}\otimes H^0(C_j,\mathbb{C})\otimes H^0(C_k,\mathbb{C}) \label{eq: decmp H1}\\
			H^2(X,\mathbb{C})&=W_2\bigoplus\bigg(\bigoplus_{\substack{1\leq i<j\leq3,\\ k\neq i,j}}\bigoplus\limits_{\chi\in I_i^0\cap I_j^0}H^1(C_i,\mathbb{C})^{\chi}\otimes H^1(C_j,\mathbb{C})^{\bar{\chi}}\otimes H^0(C_k,\mathbb{C})\bigg)\label{eq: decmp H2}\\
			H^3(X,\mathbb{C})&=W'_3\bigoplus\bigg(\bigoplus_{\substack{\chi_i\in I_i^0,1\leq i\leq3,\\ \chi_1\chi_2\chi_3=\chi_0}}H^1(C_1,\mathbb{C})^{\chi_1}\otimes H^1(C_2,\mathbb{C})^{\chi_2}\otimes H^1(C_3,\mathbb{C})^{\chi_3}\bigg)\label{eq: decmp H3}
		\end{align}
		where $W'_3=\bigoplus_{1\leq i,j\leq3}H^2(C_i,\mathbb{C})\otimes H^1(C_j,\mathbb{C})^{1_G}\otimes H^0(C_k,\mathbb{C})$
		here $I_i^0=I_i\cup\{1_G\}$ for each $1\leq i\leq 3$. From formula (\ref{eq: decmp H1}), $\mathrm{Lf}(X)=\Ker(\rho_1)$. By formula (\ref{eq: decomposition of V otimes W}),
		\begin{align*}&H^2(X,\CC)^{\chi\times\bar{\chi}\times 1_G}=H^1(C_1,\CC)^{\chi}\ot H^1(C_2,\CC)^{\bar{\chi}}\ot H^0(C_3,\CC)\\
			&H^3(X,\CC)^{\chi\times\bar{\chi}\times 1_G}=H^1(C_1,\CC)^{\chi}\ot H^1(C_2,\CC)^{\bar{\chi}}\ot H^1(C_3,\CC)^{1_G}\\
			&H^3(X,\CC)^{\chi_1\times\chi_2\times\chi_3}=H^1(C_1,\CC)^{\chi_1}\ot H^1(C_2,\CC)^{\chi_2}\ot H^1(C_3,\CC)^{\chi_3}
		\end{align*}
		where $\chi\times\bar{\chi}\times 1_G\in A_2$  and $\chi_1\times\chi_2\times\chi_3\in A_3$. Therefore, $$H^i(X,\mathbb{C})^{\psi}\neq 0\text{ if and only if }\psi\in A.$$ Compare formula (\ref{eq: decmp H2}) and formula (\ref{eq: decmp H3}), we see that $\Ker(\rho_3)\subset\Ker(\rho_2)$. By Proposition \ref{prop: AutQ to LfX}, $\mathrm{Aut}_{\mathbb{Q}}(X)=\Ker(\rho_3)=(\cap_{\psi\in A}\Ker(\psi))/K\Delta_G.$
		\qed\end{proof}
	
	\subsection{Translate to a group-theoretic problem}\label{sec:5.4}
	By Lemma \ref{lem: representation of LfX}, the numerically trivial automorphism group $\mathrm{Aut}_{\mathbb{Q}}(X)$ can be computed from the algebraic data $\mathcal{A}=(G,K_1,K_2,K_3,\cV_1,\cV_2,\cV_3)$. The rest problem is determining the set $A$ of admissible characters for $\mathcal{A}$. Now we fix some $\sigma_i\in\Sigma_i$ for each $1\leq i\leq 3$, and consider the $7$-tuple $(G,K_1,K_2,K_3,\sigma_1,\sigma_2,\sigma_3)$. In general we can consider the following datum.
	\begin{definition}\label{def: quasi algebraic data}
		Let $G$ be a finite abelian group with the identity $1$. Given three cyclic subgroups $K_1,K_2,K_3$ of $G$ and three nontrivial elements $\sigma_1,\sigma_2,\sigma_3$ in $G$. A 7-tuple $(G,K_1,K_2,K_3,\sigma_1,\sigma_2,\sigma_3)$ is called a \emph{qausi algebraic data} if it satisfies the following conditions:
		\begin{enumerate}[(1)]
			\item $K_i\cap K_j = \{1\}$ for each $1\leq i<j\leq3$;
			\item $\sigma_i\notin K_i$ for all $1\leq i\leq 3$;
			\item $\langle\sigma_1K_1\rangle\cap\langle\sigma_2K_2\rangle\cap\langle\sigma_3K_3\rangle=\{1\}$
		\end{enumerate}
		The condition (3) will be called the \emph{freeness condition}.\\
		\indent  We say that a quasi algebraic data $(G,K_1,K_2,K_3,\sigma_1,\sigma_2,\sigma_3)$ is induced from an algebraic data $\mathcal{A}=(G,K_1,K_2,K_3,\cV_1,\cV_2,\cV_3)$ if $\sigma_i\in\Sigma_i$ for all $1\leq i\leq 3$. \\
		\indent We define the admissible set relative to $\sigma_1,\sigma_2,\sigma_3$ as following $$A(\sigma_1,\sigma_2,\sigma_3):=\{\chi_1\times\chi_2\times\chi_3\in A| \chi_i(\sigma_i)\neq 1 \text{ if } \chi_i\neq 1_G, 1\leq i\leq 3\}$$ 
	\end{definition}
	Since $K_i\subseteq\mathrm{Ker}(\chi_i)$ for all $1\leq i\leq 3$, $\chi_i$ is a character of the quotient group $G/K_i$ for all $1\leq i\leq 3$. Hence the set $A(\sigma_1,\sigma_2,\sigma_3)$ depends only on the cosets $\sigma_1K_1,\sigma_2K_2,\sigma_3K_3$. For this reason, we define a equivalence between the set of quasi algebraic datum induced from a fixed algebraic data $\mathcal{A}$: 
	$$(G,K_1,K_2,K_3,\sigma_1,\sigma_2,\sigma_3)\equiv (G,K_1,K_2,K_3,\sigma_1',\sigma_2',\sigma_3')\iff \sigma_iK_i=\sigma_i'K_i \text{ for all }i$$
	So equivalent quasi algebraic datum have the same admissible set $A(\sigma_1,\sigma_2,\sigma_3)$. 
	
	Set $G(\sigma_1,\sigma_2,\sigma_3):=\cap_{\psi\in A(\sigma_1,\sigma_2,\sigma_3)}\mathrm{Ker}(\psi)$. Since for each $\chi_1\times\chi_2\times\chi_3\in A(\sigma_1,\sigma_2,\sigma_3)$, $K\Delta_G\subseteq\Ker(\chi_1\times\chi_2\times\chi_3)$ and thus $K\Delta_G\leq G(\sigma_1,\sigma_2,\sigma_3)$. By Proposition \ref{prop: AutQ to LfX}, there is an injective homomorphism $$\mathrm{Aut}_{\mathbb{Q}}(X)\hookrightarrow G(\sigma_1,\sigma_2,\sigma_3)/K\Delta_G.$$ Theorem \ref{thm: threefold isogenous AutQ}, (3) follows  from the following result.
	\begin{theorem}\label{thm: G/KD is 2 elementary abelian group}
		Under the above notations. Then the group $G(\sigma_1,\sigma_2,\sigma_3)/K\Delta_G$  is a $2$-elementary abelian group.\end{theorem}
	
	\begin{proof} Given an element $(\tau_1,\tau_2,\tau_3)\in G(\sigma_1,\sigma_2,\sigma_3)$, let $d$ be the smallest integer such that $(\tau_1^d,\tau_2^d,\tau_3^d)\in K\Delta_G$. We need to show that $d\leq 2$.
		
		Set $\tau_1'=\tau_1\tau_3^{-1},\tau_2'=\tau_2\tau_3^{-1}$, we have $(\tau_1',\tau_2',1)K\Delta_G=(\tau_1,\tau_2,\tau_3)K\Delta_G$. Observe that if $(g_1,g_2,1)\in K\Delta_G$, then $g_1\in K_1K_3$ and $g_2\in K_2K_3$.  Let $d_{i3}$ be the smallest integer such that $\tau_i'^{d_{i3}}\in K_iK_3$, then it is easy to see that $d=[d_{13},d_{23}]$ for $i=1,2$. Thus there are three integers $d_{12},d_{13},d_{23}$ such that any two of them have smallest common multiple $d$, i.e., $d_{ij}$ is the smallest integer such that $(\tau_i\tau_j^{-1})^{d_{ij}}\in K_iK_j$ for each $1\leq i<j\leq 3$.
		
		If $d_{13}\geq 3$, then $\tau_1'\notin K_1K_3$. Consider the set of characters $$I=\{\chi\in (G/K_1K_3)^*|\chi(\tau_1')\neq 1\}.$$ Applying the property $P_{1,3}$ of Lemma \ref{lem: key lemma property of P i j}, $$\chi_1\chi\times\chi_2\times\chi_3\overline{\chi}\in A(\sigma_1,\sigma_2,\sigma_3)$$
		for some $\chi\in I$ and $\chi_1\times \chi_2\times\chi_3\in A(\sigma_1,\sigma_2,\sigma_3).$ However, this shows that $\chi_1(\tau_1')\chi_2(\tau_2')=\chi_1(\tau_1')\chi(\tau_1')\chi_2(\tau_2')=1$, so we get $\chi(\tau_1')=1$, which contradicts to $\chi(\tau_1')\neq 1$. So we conclude that $d_{13}\leq 2$. For the same reason, we can prove that $d_{13},d_{23}\leq 2$. Thus we get $d\leq 2$.
		\qed\end{proof} 
	
	\begin{corollary}\label{cor: class III}
		If $X$ belongs to class III, then $\mathrm{Aut}_{\mathbb{Q}}(X)\cong (\mathbb{Z}_2)^k$ for $k=0,1,2$.
	\end{corollary}
	\begin{proof}
		By Remark \ref{rem: bound for threefold isogenous}, the order of $\mathrm{Aut}_{\mathbb{Q}}(X)$ is at most $4$, and by Theorem \ref{thm: G/KD is 2 elementary abelian group}, $\mathrm{Aut}_{\mathbb{Q}}(X)\cong (\mathbb{Z}_2)^k$ for $k=0,1,2$.
		\qed\end{proof}
	
	\subsubsection{Configurations of qausi algebraic datum}\label{sec:5.5}Given a quasi algebaric data $\mathcal{A}$. To construt admissible characters in $A(\sigma_1,\sigma_2,\sigma_3)$, we need to inverstigate the poset structure of subgroups $K_i,K_jK_k, H:=K_1K_2K_3$ of $G$ and the incidence relation between $\sigma_i$ and these subgroups.
	\begin{definition}\label{def: configuration}
		Let $\mathcal{P}(\mathcal{A}):=\{\{1\},K_i,K_jK_k,H|1\leq i\leq3,1\leq j<k\leq 3\}$ be a set of subgroups of $G$, the partial oder on $\mathcal{P}(\mathcal{A})$ is given by $U\leq V$ iff $U\subseteq V$ for $U,V\in\mathcal{P}(\mathcal{A})$. Set $\mathcal{P}(\mathcal{A})^*:=\mathcal{P}(\mathcal{A})\cup\{G\}$ with the partial order $U\leq G$ for all $U\in \mathcal{P}(\mathcal{A})$. We say that $\mathcal{A}$ is of \emph{general type} if all the subgroups $K_i,K_jK_k,H$ are distinct, otherwise it is called \emph{special type}. The \emph{configuration} associated to the algebraic data $\mathcal{A}$ is a map:
		\begin{equation*}C^{\mathcal{A}}\co \{1,2,3\}\rightarrow \mathcal{P}(\mathcal{A})^*
		\end{equation*}
		defined by $C^{\mathcal{A}}(i)$ being the smallest element $U\in\mathcal{P}(\mathcal{A})^*$ containing $\sigma_i$.
		
		If $\mathcal{A}$ is of general type, we call $\mathcal{A}$ \emph{reduced} if  $$C^{\mathcal{A}}(i)\in\{K_j,K_k,K_jK_k,G\}$$ for $\{i,j,k\}=\{1,2,3\}$.
		
		If $\mathcal{A}$ is of special type with $H=K_iK_j$, we call $\mathcal{A}$ \emph{reduced} if $$C^{\mathcal{A}}(i)\in\{K_j,G\},C^{\mathcal{A}}(j)\in\{K_i,G\}\text{ and }C^{\mathcal{A}}(k)\in\{K_i,K_j,H,G\}$$
		for $\{i,j,k\}=\{1,2,3\}$.
	\end{definition}
	
	We can visualize a poset as a graph whose vertexes are subgroups, two vertexes adjoint one edge iff the corresponding subgroups have strict order relation, and we draw this graph from top to bottum with vertexes corresponding to subgroups from $H$ to $\{1\}$. The following poset $\mathcal{P}(\mathcal{A})$ is the bigest one when $\mathcal{A}$ is of general type, i.e., all subgroups $K_i,K_jK_k, H$ are distinct.
	
	\begin{center}
		\begin{tikzpicture}
			\filldraw (0,0) circle (2pt) node[right]{$K_1$} -- (-0.25,1)circle (2pt) node[left]{$K_1K_2$}-- (1,2)circle (2pt) node[above]{$H$}
			(-0.25,1) --(1,0) circle (2pt) node[right]{$K_2$}--(2.25,1) node[right]{$K_2K_3$}(0,0)--(1,1)circle  (2pt) node[right]{$K_1K_3$}--(2,0) (1,1)--(1,2)
			(1,2)--(2.25,1) circle (2pt)--(2,0)circle (2pt) node[right]{$K_3$}
			(0,0)--(1,-0.5) circle(2pt) node[below]{$\{1\}$} (1,0)--(1,-0.5) (2,0)--(1,-0.5);
		\end{tikzpicture}
	\end{center}  
	
	\indent If $\mathcal{A}$ is of special type with $H=K_1K_2$, then $H$ is isomorphic to $K_1\times K_2$, and $K_3$ is a diagonal subgroup of $H$, i.e., $K_1\cap K_3=\{1\}, K_2\cap K_3=\{1\}$. There are three cases of $\mathcal{A}$, they are:\begin{enumerate}[(1)]
		\item  $K_1K_3,K_2K_3$ are both proper subgroups of $H$;
		\item $H=K_2K_3$ and $K_1K_3$ is a proper subgroup of $H$;
		\item $H=K_1K_3=K_2K_3$.
	\end{enumerate} 
	We draw the corresponding graphes as bellow respectively.
	\begin{center}
		\begin{tikzpicture}
			\filldraw (0,0) circle (2pt) node[right]{$K_1$}--(0.5,0.5)circle (2pt)node[left]{$K_1K_3$}--(1,1)circle (2pt)node[above]{$H$}
			(0.5,0.5)--(1,0)circle (2pt)node[right]{$K_3$}--(1.5,0.5)circle (2pt)node[right]{$K_2K_3$}--(2,0)circle (2pt)node[right]{$K_2$}
			(1,1)--(1,0) (1,1)--(1.5,0.5) 
			(0,0)--(1,-0.5) circle(2pt) node[below]{$\{1\}$} (1,0)--(1,-0.5) (2,0)--(1,-0.5);
		\end{tikzpicture}
		\begin{tikzpicture}
			\filldraw (0,0) circle (2pt) node[right]{$K_1$}--(0.5,0.5)circle (2pt)node[left]{$K_1K_3$}--(1,1)circle (2pt)node[above]{$H$}
			(0.5,0.5)--(1,0)circle (2pt)node[right]{$K_3$}
			(1,1)--(1,0) (1,1)--(2,0) circle (2pt)node[right]{$K_2$}
			(0,0)--(1,-0.5) circle(2pt) node[below]{$\{1\}$} (1,0)--(1,-0.5) (2,0)--(1,-0.5);
		\end{tikzpicture}
		\begin{tikzpicture}
			\filldraw (0,0) circle (2pt) node[right]{$K_1$}--(1,1)circle (2pt)node[above]{$H$} (1,1)--(1,0)circle (2pt)node[right]{$K_3$}
			(1,1)--(2,0) circle (2pt)node[right]{$K_2$}
			(0,0)--(1,-0.5) circle(2pt) node[below]{$\{1\}$} (1,0)--(1,-0.5) (2,0)--(1,-0.5);
		\end{tikzpicture}
	\end{center} 
	For each class of quasi algebraic datum, we want to choose a simple representative which is a reduced quasi algebraic data, we have the following proposition.
	\begin{proposition}\label{prop: reduce to reduced algebraic data}
		Each algebraic data $\mathcal{A}$ is equivalent to a reduced algebraic data $\mathcal{A}'$.
	\end{proposition}
	\begin{proof} Suppose $\mathcal{A}$ is of general type and not reduced, then by Definition \ref{def: configuration}, there is some $i$ such that $C(i)^{\mathcal{A}}=K_iK_j$ or $C(i)^{\mathcal{A}}=H$, it follows that $\sigma_i\in K_iK_j$ or $\sigma_i\in H$. If $\sigma_i\in K_iK_j$, we can write $\sigma_i=k_i\sigma_i'$ for some $k_i\in K_i$ and $\sigma_i'\in K_j$, hence $\sigma_iK_i=\sigma_i'K_i$, so we can replace $\sigma_i$ in $\mathcal{A}$ by $\sigma_i'$ and get a new algebraic data $\mathcal{A}'$ such that $C^{\mathcal{A}'}(i)=K_j$. Similarly for the case $C^{\mathcal{A}}(i)=H$, we will get a new algebraic data $\mathcal{A}'$ with $C^{\mathcal{A}'}(i)=K_jK_k$, where $\{i,j,k\}=\{1,2,3\}$. Continoue this procedure, we will get a reduced algebraic data $\mathcal{A}'$ which is equivalent to the original one. The same argument works for $\mathcal{A}$ being of special type.
		\qed\end{proof}
	
	\subsubsection{Construction of admissible characters}\label{sec:5.6}
	\begin{lemma}\label{lem: key lemma property of P i j}
		Let $\mathcal{A}$ be a quasi algebraic data, then the admissible set $A(\sigma_1,\sigma_2,\sigma_3)$ is nonempty.
		
		Moreover, the set $A(\sigma_1,\sigma_2,\sigma_3)$ satisfies the property $P_{i,j}$: for $\{i,j,k\}=\{1,2,3\}$ and $g\notin K_iK_j$ with $o(gK_iK_j)\geq 3$, $$\chi_i\chi\times\chi_j\overline{\chi}\times\chi_k\in A(\sigma_1,\sigma_2,\sigma_3)$$
		for some $\chi_i\times\chi_j\times\chi_k\in A(\sigma_1,\sigma_2,\sigma_3)$ and $\chi\in(G/K_iK_j)^*$ such that $\chi(g)\neq 1$.
	\end{lemma}
	
	\begin{proof}
		\indent By Proposition \ref{prop: reduce to reduced algebraic data}, we can assume $\mathcal{A}$ is reduced. If $\mathcal{A}$ is of general type. Recall that $C^{\mathcal{A}}(i)$ takes values in the set $\{K_j,K_k,K_jK_k,G\}$. The proof will be divided into three steps according to the number $\nu$ of $C^{\mathcal{A}}(i)=G$.
		\paragraph{General type $\nu = 2$.} Suppose $C^{\mathcal{A}}(1)=C^{\mathcal{A}}(3)=G$, then $\sigma_1\notin K_1K_3$ and  $\sigma_3\notin K_1K_3$. There is a character  $\phi\in (G/K_1K_3)^*$ such that \begin{equation*}
			\phi(\sigma_1)\neq 1\text{ and }\phi(\sigma_3)\neq 1.
		\end{equation*} Thus $\phi\times 1_G\times\overline{\phi}\in A(\sigma_1,\sigma_2,\sigma_3)$, which is of weight $2$.\\
		\indent  Moreover, if $C^{\mathcal{A}}(i)\neq G$ for some $1\leq i\leq 3$. Suppose that $$C^{\mathcal{A}}(1)=K_3\text{ and } C^{\mathcal{A}}(2)=C^{\mathcal{A}}(3)=G.$$ We have
		\begin{align*}
			&\beta_3(\sigma_1)\neq1\text{ for some }\beta_3\in(H/K_1K_2)^*\\
			&\phi(\sigma_2)\neq1\text{ and }\phi(\sigma_3)\neq 1\text{ for some }\phi\in(G/K_2K_3)^*
		\end{align*}
		Note that $\beta_3^G(\sigma_2)=0$ provided that $\sigma_2\notin H$, by Proposition \ref{prop: consituent of induced character such that c chi neq 1}, there is a consituent $\psi_3$ of $\beta_3^G$ such that $\phi(\sigma_2)\psi_3(\sigma_2)\neq1$. It follows that $$\psi_3\times\overline{\phi\psi_3}\times\phi\in A(\sigma_1,\sigma_2,\sigma_3),$$ which is of weight $3$.
		\paragraph{General type $\nu=1$.}Suppose $C^{\mathcal{A}}(3)=G$. Since $$C^{\mathcal{A}}(1)\in\{K_2,K_3,K_2K_3\}\text{ and } C^{\mathcal{A}}(2)\in\{K_1,K_3,K_1K_3\},$$ there are nine possible choices of $C^{\mathcal{A}}(1)$ and $C^{\mathcal{A}}(2)$. In each choice, the constructions are similar, we just illustrate one case. Suppose $$C^{\mathcal{A}}(1)=K_3,C^{\mathcal{A}}(2)=K_1\text{ and }C^{\mathcal{A}}(3)=G.$$ Since $\sigma_3\notin H$, there is a character $\phi\in(G/H)^*$ such that $\phi(\sigma_3)\neq1$. We first take a character $\beta_1\in(H/K_2K_3)^*$ with $\beta_1(\sigma_2)\neq 1$, as $\sigma_3\notin H$, we have $\beta_1^G(\sigma_3)=0$, by Proposition \ref{prop: consituent of induced character such that c chi neq 1}, there is a constituent $\psi_1$ of $\beta_1^G$ such that 
		\begin{equation}\label{eq: phi psi neq 1}
			\phi(\sigma_3)\psi_1(\sigma_3)\neq1.
		\end{equation}
		Take a character $\psi_3\in(G/K_1K_2)^*$ such that $\psi_3(\sigma_1)\neq 1$, then we obtain
		\begin{equation}\label{eq: phi psi1 psi3=beta neq 1}
			\phi(\sigma_2)\psi_1(\sigma_2)\psi_3(\sigma_2)=\beta_1(\sigma_2)\neq1,\end{equation}
		as $\sigma_2\in K_1$ and $K_1\subseteq\Ker(\phi),K_1\subseteq\Ker(\psi_3)$. Combine (\ref{eq: phi psi neq 1}) and (\ref{eq: phi psi1 psi3=beta neq 1}), we deduce that $\psi_3\times\overline{\phi\psi_1\psi_3}\times\phi\psi_1\in A(\sigma_1,\sigma_2,\sigma_3)$.
		\paragraph{General type $\nu=0$.}We divide three cases.
		\paragraph{Case (a):} Suppose that $$C^{\mathcal{A}}(1)\in\{K_3,K_2K_3\},C^{\mathcal{A}}(2)\in\{K_1,K_1K_3\}\text{ and }C^{\mathcal{A}}(3)\in\{K_2,K_1K_2\}.$$ We can take characters $\psi_1\in(G/K_2K_3)^*,\psi_2\in(G/K_1K_3)^*,\psi_3\in(G/K_1K_2)^*$ such that $\psi_i(\sigma_i)\neq1$ for all $i$. It is easy to verify that 
		\begin{equation*}
			\psi_2\psi_3\times\psi_1\overline{\psi_3}\times\overline{\psi_1\psi_2}\in A(\sigma_1,\sigma_2,\sigma_3)
		\end{equation*}
		\paragraph{Case (b):} Suppose that $$C^{\mathcal{A}}(1)\in\{K_3,K_2K_3\},C^{\mathcal{A}}(2)\in\{K_3,K_1K_3\}\text{ and }C^{\mathcal{A}}(3)\in\{K_1,K_1K_2\}.$$ We can take characters $\psi_1\in(G/K_2K_3)^*,\psi_3\in(G/K_1K_2)^*$ such that $$\psi_1(\sigma_3)\neq1,\psi_3(\sigma_1)\neq1\text{ and }\psi_3(\sigma_2)\neq1.$$ Therefore 
		\begin{equation*}
			\psi_3\times\psi_1\overline{\psi_3}\times\overline{\psi_1}\in A(\sigma_1,\sigma_2,\sigma_3).
		\end{equation*}
		\paragraph{Case (c):} Suppose that $$C^{\mathcal{A}}(1)=K_2,C^{\mathcal{A}}(2)=K_1\text{ and }C^{\mathcal{A}}(3)=K_1K_2.$$
		Then we can write $\sigma_3=k_1k_2$ for some $k_i\in K_i$ for $i=1,2$. Set $m=o(\sigma_1),n=o(\sigma_2),m'=o(k_1)\text{ and }n'=o(k_2)$.\\ 
		\indent Since $K_2K_3$ and $K_1K_3$ are proper subgroups of $H$, there are characters $\beta_1\in (H/K_2K_3)^*$ and $\beta_2\in(H/K_1K_3)^*$ whose restriciton on $K_1$ and $K_2$ are primitive characters, respectively. If $\max\{m,n,m',n'\}\geq 3$, since $K_2\subseteq\Ker(\beta_1)$ and $K_1\subseteq\Ker(\beta_2)$, we can apply Lemma \ref{lem: if max geq 3 construction of adm characters} to the group $K_1K_2$ with $g_1=\sigma_1,g_2=\sigma_2,h=\sigma_3,\alpha_1=(\beta_2)_{K_1K_2},\alpha_2=(\beta_1)_{K_1K_2}$, then we have 
		$$\begin{cases}
			\beta_2^s(\sigma_1)\neq 1\text{ and }\beta_1(\sigma_3)\beta_2^s(\sigma_3)\neq 1 &\text{ if }m\geq3\text{ or }m'\geq3,\\
			\beta_1^s(\sigma_2)\neq 1\text{ and }\beta_1^s(\sigma_3)\beta_2(\sigma_3)\neq 1 &\text{ if }n\geq3\text{ or }n'\geq3,
		\end{cases}$$ for some $s=1,2,3,5$. It follows that
		\begin{equation*}
			\begin{cases}
				\psi_2^s\times\psi_1\times\overline{\psi_1\psi_2^s}\in A(\sigma_1,\sigma_2,\sigma_3) & \text{ if }m\geq3\text{ or }m'\geq3,\\
				\psi_2\times\psi_1^s\times\overline{\psi_1^s\psi_2}\in A(\sigma_1,\sigma_2,\sigma_3) & \text{ if }n\geq3\text{ or }n'\geq3,
			\end{cases}
		\end{equation*} for some $s=1,2,3,5$, where $\psi_1,\psi_3$ are characters of $G$ whose restricition on $H$ are equal to $\beta_1,\beta_3$, respectively. Now we can assume that $m=n=m'=n'=2$, then we can write $\sigma_3 = \sigma_1\sigma_2$. We can verify that $\sigma_3\in \langle\sigma_1 K_1\rangle\cap\langle\sigma_2 K_2\rangle\cap\langle \sigma_3 K_3\rangle$. For this case the corresponding quasi algebraic data violates the freeness condition of Definition \ref{def: quasi algebraic data}.
		\paragraph{Special type.}Now we assume that $\mathcal{A}$ is of special type $H=K_1K_2$. Recall that $$C^{\mathcal{A}}(1)\in\{K_2,G\},C^{\mathcal{A}}(2)\in\{K_1,G\}\text{ and }C^{\mathcal{A}}(3)\in\{K_1,K_2,H,G\}.$$ If $\nu=1,2$, we can construct admissible character as for the general type $\nu=1,2$. For the case $C^{\mathcal{A}}(1)=K_2,C^{\mathcal{A}}(2)=K_1,C^{\mathcal{A}}(3)\in\{K_1,K_2\}$, we can construct admissible character as the case (b) of $\mathcal{A}$ being of general type with $\nu=0$. So we can assume that $$C^{\mathcal{A}}(1)=K_2,C^{\mathcal{A}}(2)=K_1\text{ and }C^{\mathcal{A}}(3)=H.$$ Let $m_1=o(\sigma_1),m_2=o(\sigma_2)$.  Sinces $H\cong K_1\times K_2$, and $K_3$ is a subgroup such that $K_1\cap K_3=\{1\}$ and $K_2\cap K_3=\{1\}$. We can apply Lemma \ref{lem: construction of bases characters}, there are characters $\alpha_i\in H^*$ for $i=1,2$ such that $$K_1\subseteq\Ker(\alpha_2),K_2\subseteq\Ker(\alpha_1)\text{ and }K_3\subseteq\Ker(\alpha_1\alpha_2),$$ and the restriction of $\alpha_i$ on $K_i$ is primitive for each $i=1,2$. Since $\sigma_3\in H$, we can write $\sigma_3=k_1k_2$ for $k_i\in K_i, i=1,2$. Similary as the case (c) of $\mathcal{A}$ being of general type with $\mu=0$, if $\max\{m_1,m_2,o(k_1),o(k_2)\}\geq 3$, then we can apply Lemma \ref{lem: if max geq 3 construction of adm characters} to the group $H$ with $g_1=\sigma_2,g_2=\sigma_1,h=\sigma_3$ and $\alpha_1,\alpha_2$, we have that either $\psi_2\times\psi^s_1\times\overline{\psi^s_1\psi_2}\in A(\sigma_1,\sigma_2,\sigma_3)$ or $\psi^s_2\times\psi_1\times\overline{\psi_1\psi^s_2}\in A(\sigma_1,\sigma_2,\sigma_3)$ for some integer $s$, where $\psi_1,\psi_2$ are any characters of $G$ whose restriciton on $H$ is equal to $\alpha_1,\alpha_2$, respectively. In the case $m_1=m_2=o(k_1)=o(k_2)$, the corresponding quasi algebraic data violates the freeness condition of Definition \ref{def: quasi algebraic data}.\\
		\paragraph{Proof of property $P_{1,3}$.} If $\sigma_1\notin K_1K_3$ and $\sigma_3\notin K_1K_3$, then we can find a character $\phi\times1_G\times\overline{\phi}\in A(\sigma_1,\sigma_2,\sigma_3)$ such that
		\begin{equation*}
			\phi\in(G/K_1K_3)^*, \phi(\sigma_1)\neq 1 \text{ and } \phi(\sigma_3)\neq 1.
		\end{equation*} Since $o(gK_1K_3)\geq 3$, we can apply lemma \ref{lem: if o geq 3 construction of adm characters} to the quotient group $G/K_1K_3$ and its three nontrivial elements $\sigma_1K_1K_3,\sigma_3K_1K_3, gK_1K_3$, there are characters $\phi,\chi\in(G/K_1K_3)^*$ such that \begin{equation*}
			\phi(\sigma_1)\neq1,\phi(\sigma_3)\neq 1,\chi(g)\neq 1,\phi(\sigma_1)\chi(\sigma_1)\neq 1\text{ and } \phi(\sigma_3)\chi(\sigma_3)\neq 1.
		\end{equation*} It follows that $\phi\times1_G\times\overline{\phi}\in A(\sigma_1,\sigma_2,\sigma_3),\chi\in I$ and $\phi\chi\times1_G\times\overline{\phi\chi}\in A(\sigma_1,\sigma_2,\sigma_3)$.\\
		\indent Since $\sigma_1\notin K_1K_3$ is equivalent to $C^{\mathcal{A}}(1)\neq K_3$, next step we consider the case $C^{\mathcal{A}}(1)= K_3$. By the constructions as above, there exists a character $\chi_1\times\chi_2
		\times\chi_3$ of weight $3$, i.e., $\chi_i\neq 1_G$ for all $1\leq i\leq 3$. Since $\sigma_1\in K_3$, for any $\chi\in I$, we have $\chi(\sigma_1)=1$, hence $\chi_1(\sigma_1)\chi(\sigma_1)=\chi_1(\sigma_1)\neq 1$. Since $o(gK_1K_3)\geq 3$, we can take $\chi\in I$ such that $\chi(g)$ is a primitive $o(gK_1K_3)$-th root of unit, thus we have $\chi^2\in I$. Observe that if $\chi_3(\sigma_3)\chi(\sigma_3)=\chi_3(\sigma_3)\chi^2(\sigma_3)=1$, then $\chi(\sigma_3)=\chi_3(\sigma_3)=1$, which contradict to $\chi_3(\sigma_3)\neq 1$. So we can choose $\chi\in I$ such that $\chi_3(\sigma_3)\chi^s(\sigma_3)\neq 1$ for some $s=1,2$. It follows that $\chi_1\chi^s\times\chi_2
		\times\chi_3\overline{\chi^s}\in A(\sigma_1,\sigma_2,\sigma_3)$ for some $s=1,2$.
		\qed\end{proof}
	
	\section{Examples}\label{sec:6}
	
	\indent  Let $n_1,n_2,n_3$ be three positive integers. In Examples \ref{ex: Z2 Z2}, \ref{ex: Z2}  and \ref{ex:Z2n}, let $G=\langle e_1,e_2,e_3\rangle\cong\ZZ_{2n_1}\times\ZZ_{2n_2}\times\ZZ_{2n_3}$ where $e_1,e_2,e_3$ are generators of $G$ of orders $2n_1,2n_2,2n_3$ respectively. Let $\phi_i$ be the characters on $G$ such that $\phi_i(e_j)=e^{2\pi\delta_{ij}\sqrt{-1}/2n_i}$, where $\delta_{ij}$ is the Kronecker symbol. For each $1\leq i\leq 3$, let $K_i=\langle e_i\rangle$. We take the following generating vector of $G/K_i$
	\begin{align*}
		\cV_1&=(e_2K_1, e_3K_1;\underbrace{\sigma_1K_1,\dots,\sigma_1K_1}_{m_1}),\\
		\cV_2&=( e_1K_2, e_3K_2;\underbrace{\sigma_2K_2,\dots,\sigma_2K_2}_{m_2}),\\
		\cV_3&=( e_1K_3, e_2K_3;\underbrace{\sigma_3K_3,\dots,\sigma_3K_3}_{m_3}).
	\end{align*}
	Where $m_i$ is an integer which is divisible by $o(\sigma_iK_i)$ for each $1\leq i\leq 3$. So the type of $\cV_i$ is $$[1;\underbrace{o(\sigma_iK_i),\dots,o(\sigma_iK_i)}_{m_i}].$$ By Riemann's existence theorem, there is an algebraic curve $C_i$ with a $G$-action $\psi_i\co G\to\Aut(C_i)$ such that $K_i=\Ker(\psi_i)$ whose quotient $C_i/G$ is an elliptic curve for each $1\leq i\leq3$. Since the set of nontrivial stabilizers of the $G$-action on the $C_i$ is $\Sigma_i=\sigma_iK_i$.
	
	\begin{example}[$\mathrm{Aut}_{\mathbb{Q}}(X)\cong \ZZ_2\times \ZZ_2$]\label{ex: Z2 Z2}\quad\\
		\indent	Take $\sigma_1=e_3^{n_3},\sigma_2=e_1^{n_1}\text{ and }\sigma_3=e_2^{n_2}$. Thus $\cap_i\Sigma_i=\{1_G\}$, and the $7$-tuple $(G,K_1,K_2,K_3,\cV_1,\cV_2,\cV_3)$ forms an algebraic datum.  By Remark \ref{rem: construction of threefold isogenous}, it  determines a threefold $X=(\curProd)/G$ isogenous to an unmixed product of curves with $q(X)=3$. It is easy to see that the corresponding admissible sets are
		$$A_2=\varnothing \text{ and } A_3=\{\phi_2^{k_2}\phi_3^{k_3}\times\phi_1^{k_1}\phi_3^{-k_3}\times\phi_1^{-k_1}\phi_2^{-k_2}|k_1,k_2,k_3\text{ are odd integers}\}.$$
		Let $\tau=(e^{x_{1}}_2e^{x_{2}}_3,e_1^{x_{3}}e_3^{x_{4}},1)\in G\times G\times G$ be a representative of an element of $\mathrm{Aut}_{\mathbb{Q}}(X)$ where $x_{i}\in \ZZ$. Then $\tau K\Delta_G\in\mathrm{Aut}_{\mathbb{Q}}(X)$ if and only if the following equations hold:
		\begin{align*}
			\phi_1(e_1)^{k_1x_{3}}\phi_2(e_2)^{k_2x_{1}}\phi_3(e_3)^{k_3(x_{2}-x_{4})}=1
		\end{align*}
		for all odd integers $k_i$. In particular from the two equations 
		\begin{align*}
			\phi_1(e_1)^{x_{3}}\phi_2(e_2)^{x_{1}}\phi_3(e_3)^{x_{2}-x_{4}}=\phi_1(e_1)^{3x_{3}}\phi_2(e_2)^{x_{1}}\phi_3(e_3)^{x_{2}-x_{4}}=1\end{align*} 
		we get $\phi_1(e_1)^{2x_{3}}=1$. Since $\phi_1(e_1)$ is a primitive $2n_1$-th root of unit, $2n_1|2x_{3}$, thus the possible value of $x_{3}$ is $n_1$ or $0$. Applying this argument again we have $x_{1}=n_2,0$ and $x_2-x_4=n_3, 0$. By Lemma \ref{lem: representation of LfX} we have $$\mathrm{Aut}_{\mathbb{Q}}(X)=\langle(e_2^{n_2},e_1^{n_1},1)K\Delta_G,(e_2^{n_2},e_3^{n_3},1)K\Delta_G\rangle\cong \ZZ_2\times \ZZ_2.$$
	\end{example}
	
	\begin{example}[$\mathrm{Aut}_{\mathbb{Q}}(X)\cong \ZZ_2$]\label{ex: Z2}\quad\\
		\indent Take $\sigma_1=e_3^{n_3},\sigma_2=e_1^{n_1}\text{ and }\sigma_3=e_1^{n_1}e_2^{n_2}$,
		for the same reason as Example \ref{ex: Z2 Z2}, the $7$-tuple $(G,K_1,K_2,K_3,\cV_1,\cV_2,\cV_3)$ forms an algebraic datum. Let $X$ be the corresponding threefold. We can see  that the corresponding admissible sets are
		\begin{align*}
			&A_2=\varnothing \text{ and }\\ &A_3=\{\phi_2^{k_2}\phi_3^{k_3}\times\phi_1^{k_1}\phi_3^{-k_3}\times\phi_1^{-k_1}\phi_2^{-k_2}|k_1,k_2\text{ are odd, }k_3\text{ is even}\}.
		\end{align*} 
		By Lemma \ref{lem: representation of LfX} we obtain $$\mathrm{Aut}_{\mathbb{Q}}(X)=\langle(e_2^{n_2},1,1)K\Delta_G\rangle\cong \ZZ_2.$$\\
		\indent Take $\sigma_1=e_3^{n_3},\sigma_2=e_1^{2}\text{ and }\sigma_3=e_2^{n_2}$. Similarly,  we have \begin{align*}
			&A_2=\varnothing \text{ and }\\ &A_3=\{\phi_2^{k_2}\phi_3^{k_3}\times\phi_1^{k_1}\phi_3^{-k_3}\times\phi_1^{-k_1}\phi_2^{-k_2}|k_1\text{ is arbitrary, }k_2,k_3\text{ are odd}\}. 
		\end{align*} By Lemma \ref{lem: representation of LfX}, it follows that $$\mathrm{Aut}_{\mathbb{Q}}(X)=\langle(e_2^{n_2},e_3^{n_3},1)K\Delta_G\rangle\cong \ZZ_2.$$ 
	\end{example}
	
	\begin{example}[Product quotient with terminal singularities]\label{ex:Z2n}\quad\\
		Let $n=n_1$ and $m=n_2=n_3$. Take $\sigma_1=e_2^me_3^m,\sigma_2=e_3^m$ and $\sigma_3=e_2^m$. Since $\sigma_1 K_2= \sigma_2 K_2, \sigma_1 K_3= \sigma_3 K_3$ and $\sigma_iK_i$ has fixed point of $G/K_i$-action on $C_i$ for all $i=1,2,3$, hence $(\sigma_1,\sigma_1,\sigma_1)$ has fixed points on $\curProd$. Therefore, the corresponding $G$-action on $\curProd$ is not free.  So the quotient $X_{m,n}=(\curProd)/G$ has singularities of type $\frac{1}{2}(1,1,1)$ which is not Gorenstein. Even though, cohohomolgies $H^*(X_{m,n},\CC)$ of $X_{m,n}$ can be identified with $H^{*}(\curProd,\CC)^G$. The same arguments in Section \ref{sec:5} apply to $X_{m,n}$, we can see that Lemma \ref{lem: representation of LfX} still holds for $X_{m,n}$. The  corresponding admissible sets are
		$$A_2=\{\phi_2^{k_2}\times1_G\times\phi_2^{-k_2},\phi_3^{k_3}\times\phi_3^{-k_3}\times1_G|k_2,k_3 \text{ are odd}\}\text{ and }A_3=\varnothing.$$
		It follows that $$\mathrm{Aut}_{\mathbb{Q}}(X_{m,n})=\langle (1,e_1,1)K\Delta_G\rangle\cong\ZZ_{2n}.$$
		For $\sigma\in G$ and $1\leq i\leq3$, we denote the set of points on $C_i$ whose stabilizer is exactly $\langle\sigma\rangle K_i$ by $\mathrm{Fix}_{C_i}(\sigma)$. By \cite[Lemma 10.4]{Bre00}, we have 
		\begin{align*}
			&|\mathrm{Fix}_{C_1}(e_2^me_3^m)|=4m^2,|\mathrm{Fix}_{C_2}(e_3^m)|=|\mathrm{Fix}_{C_3}(e_2^m)|=4mn;\\
			&|\mathrm{Fix}_{\curProd}(\sigma)|=4^3m^4n^2.
		\end{align*}
		Therefore, the basket of singularities of $X_{m,n}$ is $$\mathrm{Sing}(X_{m,n})=\{\frac{4^3m^4n^2}{4m^2n}\cdot\frac{1}{2}(1,1,1)\}=\{16m^2n\cdot\frac{1}{2}(1,1,1)\}.$$ It is easy to see that $K_{X_{m,n}}.c_2(X_{m,n})=\frac{3\cdot2^3\prod_{i=1}^3(g(C_i)-1)}{|G|}= 24m^2n$. By Riemann-Roch formula for singular varieties \cite[Corollary 10.3]{Re85}, we have 
		$$\chi(\omega_{X_{m,n}})=\frac{1}{24}K_{X_{m,n}}. c_2(X_{m,n})-16m^2n\cdot\frac{1}{16}= 0.$$ According to the classification of threefolds with vanishing holomorphic Euler characteristic by J. A. Chen, O. Debarre and Z. Jiang \cite{CDJ14}, the variety $X_{m,n}$ belongs to the examples constructed by R. Lararsfeld and L. Ein \cite[Example 1.13]{EL97}.
	\end{example}
	
	\begin{example}[$\mathrm{Aut}_{\mathbb{Q}}(X)\cong\mathbb{Z}_2$ and $q(X)>3$]\label{ex: q > 3 Z2}
		Let $G=\langle e_1,e_2,e_3\rangle\cong \mathbb{Z}_2\times\mathbb{Z}_2\times\mathbb{Z}_2$, $K_1=K_2=\{1\}$ and $K_3=\langle e_1,e_2\rangle $, $\phi_i$ is the dual character of $e_i$. We can take generating vectors for $G$ as the following:
		\begin{align*}
			&\cV_1=(e_2,e_3;e_1,e_1)\\
			&\cV_2=(e_1,e_3;e_2, e_2).
		\end{align*}
		Take  $\cV_3=(4\cdot e_3K_3;e_3K_3,e_3K_3)$ a generating vector for $G/K_3$. It is easy to see that the $7$-tuple $(G,K_1,K_2,K_3,\cV_1,\cV_2,\cV_3)$ forms an algebraic datum. Then  the corresponding admissible sets are
		\begin{align*}
			&A_2=\{\phi_1\phi_2\times\phi_1\phi_2\times 1_G\}\\
			&A_3=\{\phi_1\phi_2\phi_3\times\phi_1\phi_2\times \phi_3,\phi_1\phi_2\times\phi_1\phi_2\phi_3\times \phi_3\}
		\end{align*}
		By Lemma \ref{lem: representation of LfX} and Corollary \ref{cor: class II} we know that $\mathrm{Aut}_{\mathbb{Q}}(X)\cong\langle e_1,e_2,1\rangle K\Delta_G$.
		
	\end{example}
	
	\begin{example}[Some $K_i$ is not cyclic and $\mathrm{Aut}_{\mathbb{Q}}(X)\cong\ZZ_2$]\label{ex: K not cyclic}\quad\\
		\indent	 Let $G=\ZZ^4_{2}$ with generators $e_1,e_2,e_3,e_4$, and let $$K_1=\langle e_4\rangle, K_2=\langle e_2\rangle, K_3=\langle e_1,e_3\rangle.$$ We take generating vectors as following 
		\begin{align*}
			&\cV_1=(e_2K_1,e_3K_1;e_1K_1, e_1K_1),\\
			&\cV_2=(e_1K_2,e_4K_2;e_3K_2, e_3K_2),\\
			&\cV_3=(e_2K_3, e_3K_3;e_2K_3,e_2K_3).
		\end{align*}
		It is easy to see that the $7$-tuple $(G,K_1,K_2,K_3,\cV_1,\cV_2,\cV_3)$ forms an algebraic datum.
		The  corresponding admissible sets are \begin{align*}
			&A_2=\{\phi_1\phi_3\times\phi_1\phi_3\times1_G\}\text{ and }\\
			&A_3=\{\phi_1\phi_2\phi_3\times\phi_1\phi_3\times\phi_2,\phi_1\phi_2\phi_3\times\phi_1\phi_3\phi_4\times\phi_2\phi_4\}
		\end{align*} 
		It follows that $\mathrm{Aut}_{\mathbb{Q}}(X)=\langle (e_3,e_1,1)K\Delta_G\rangle\cong\ZZ_2$.
	\end{example}
	
	\begin{question}Let $X$ be a threefold isogenous to a product of curves, not necessary unmixed type, with maximal Albanese dimension. Does $\mathrm{Aut}_{\mathbb{Q}}(X)\cong\ZZ_2^k$ for some $k=0,1,2$?
	\end{question}


	\section*{Acknowledgements}
		The author would like to thank his adivisor Jinxing Cai for constant support and encouragement, and he also would like to thank Wenfei Liu, Lei Zhang, Zhan Li and Yifan Chen for useful suggestions and discussions during the conference at Xiamen University in June, 2018. He would like to thank Xiamen University for its hospitality.

\bibliographystyle{alpha}   
	\bibliography{Ref}   
	
	%
	%
	
\end{document}